\documentclass[11pt, a4paper]{article}

\usepackage{xcolor}

\definecolor{lifex}{HTML}{f60248}
\newcommand{\lifex}{\texttt{life\textsuperscript{\color{lifex}{x}}}}
\newcommand{\lifexep}{\lifex{}\texttt{-ep}}
\newcommand{\lifexfiber}{\lifex{}\texttt{-fiber}}
\newcommand{\lifexcfd}{\lifex{}\texttt{-cfd}}
\newcommand{\zenodourl}{\url{https://doi.org/10.5281/zenodo.8085266}}

\newcommand{\dealii}{\texttt{deal.II}}

\definecolor{codegreen}{rgb}{0,0.6,0}
\definecolor{codegray}{rgb}{0.5,0.5,0.5}
\definecolor{codepurple}{rgb}{0.58,0,0.82}
\definecolor{backcolor}{rgb}{0.95,0.95,0.92}
\colorlet{linkcolor}{blue!50!black}

\usepackage[T1]{fontenc}
\usepackage[utf8]{inputenc}
\usepackage{csquotes}
\usepackage[english]{babel}
\usepackage{tipa}

\usepackage{amsmath}
\usepackage[overload]{empheq}
\usepackage{relsize}
\usepackage{siunitx}

\usepackage{authblk}

\usepackage[backend=biber, giveninits=true, maxbibnames=20, sorting=none]{biblatex}
\usepackage{hyperref}
\hypersetup{
    colorlinks,
    linktoc=page,
    linkcolor=linkcolor,
    citecolor=linkcolor,
    urlcolor=linkcolor
}
\usepackage[sort]{cleveref}
\usepackage{url}
\usepackage{acronym}

\newcommand\blankfootnote[1]{%
    \begingroup
    \renewcommand\thefootnote{}\footnote{#1}%
    \addtocounter{footnote}{-1}%
    \endgroup
}

\usepackage{listings}
\crefname{lstlisting}{listing}{listings}
\Crefname{lstlisting}{Listing}{Listings}

\lstdefinestyle{mystyle}{
    backgroundcolor=\color{backcolor},
    commentstyle=\color{codegreen},
    keywordstyle=\color{lifex},
    numberstyle=\tiny\color{codegray},
    stringstyle=\color{codepurple},
    basicstyle=\ttfamily\footnotesize,
    breakatwhitespace=false,
    breaklines=true,
    captionpos=b,
    keepspaces=true,
    numbers=none,
    numbersep=5pt,
    showspaces=false,
    showstringspaces=false,
    showtabs=false,
    tabsize=4,
    frame=shadowbox
}

\lstdefinelanguage{prm}{
    keywords={subsection, set, end},
    comment=[l]{\#}
}

\lstdefinelanguage{xml}{
    keywords={value, default\_value, documentation, pattern, pattern\_description}
}

\lstdefinelanguage{json}{
    keywords={value, default\_value, documentation, pattern, pattern\_description}
}

\DeclareUnicodeCharacter{0308}{\"}
\DeclareUnicodeCharacter{0301}{\'}

\newcommand{\fZero}{{\mathbf{f}_0}}
\newcommand{\sZero}{{\mathbf{s}_0}}
\newcommand{\nZero}{{\mathbf{n}_0}}

\newcommand{\Chim}{\chi_\mathrm{m}}
\newcommand{\Cm}{C_\mathrm{m}}
\newcommand{\Iion}[1]{{\mathcal{I}_{\mathrm{ion}}^{#1}}}
\newcommand{\Iapp}{{\mathcal{I}_{\mathrm{app}}}}
\newcommand{\DiffTens}[1]{\boldsymbol{D}^{#1}}
\newcommand{\IionDIM}[1]{{\widehat{\mathcal{I}}_{\mathrm{ion}}^{#1}}}
\newcommand{\IappDIM}{{\widehat{\mathcal{I}}_{\mathrm{app}}}}
\newcommand{\DiffTensDIM}[1]{\widehat{\boldsymbol{D}}^{#1}}
\newcommand{\Pot}{u}
\newcommand{\Gating}[1]{\boldsymbol{w}^{#1}}
\newcommand{\PotStart}[1]{u_0^{#1}}
\newcommand{\GatingStart}[1]{\boldsymbol{w}_0^{#1}}
\newcommand{\RhsGating}[1]{\boldsymbol{H}^{#1}}
\newcommand{\sigmal}[1]{\sigma_{\text{l}}^{#1}}
\newcommand{\sigmat}[1]{\sigma_{\text{t}}^{#1}}
\newcommand{\sigman}[1]{\sigma_{\text{n}}^{#1}}

\newcommand{\Domain}[1]{\Omega_{#1}}
\newcommand{\ClosedDomain}[1]{\overline{\Omega}_{#1}}
\newcommand{\Normal}[1]{\mathbf{n}_{#1}}
\newcommand{\Boundary}[1]{\Gamma_{#1}}

\newcommand{\NumIonicVariables}[1]{M_{#1}}
\newcommand{\NumDomains}{N}

\newcommand{\APot}{\mathsf u}

\newcommand{\AIion}{\mathsf s}
\newcommand{\AIapp}{\mathsf f}

\newcommand{\Mass}{\mathsf M}
\newcommand{\Stiff}{\mathsf K}

\newcommand{\alphabdf}{\alpha_\text{BDF}}
\newcommand{\bdf}[1]{_{\text{BDF},#1}}
\newcommand{\ext}[1]{_{\text{EXT},#1}}

\newcommand{\delmb}[1]{\ignorespaces}

\begin{document}
\setcounter{page}{1}

\title{\lifexep{}: a robust and efficient software for cardiac electrophysiology simulations}

\author[a,b]{Pasquale C. Africa}
\author[a,*]{Roberto Piersanti}
\author[a]{Francesco Regazzoni}
\author[a]{Michele Bucelli}
\author[a,c]{Matteo Salvador}
\author[a]{Marco Fedele}
\author[a]{Stefano Pagani}
\author[a]{Luca Dede'}
\author[a,d]{Alfio Quarteroni}

\affil[a]{\footnotesize MOX, Department of Mathematics, Politecnico di Milano, Italy}
\affil[b]{\footnotesize mathLab, Mathematics Area, SISSA International School for Advanced Studies, Trieste, Italy}
\affil[c]{\footnotesize Institute for Computational and Mathematical Engineering, Stanford University, California, USA}
\affil[d]{\footnotesize Institute of Mathematics, \'Ecole Polytechnique Fédérale de Lausanne, Switzerland (Professor emeritus)}
\affil[*]{Corresponding Author: roberto.piersanti@polimi.it}

\date{03 August 2023}

\maketitle

\blankfootnote{Official website: \url{https://lifex.gitlab.io/}}
\blankfootnote{Download link: \url{https://doi.org/10.5281/zenodo.8085266}}

\begin{abstract}
Simulating the cardiac function requires the numerical solution of multi-physics and multi-scale mathematical models. This underscores the need for streamlined, accurate, and high-performance computational tools. Despite the dedicated endeavors of various research teams, comprehensive and user-friendly software programs for cardiac simulations, capable of accurately replicating both physiological and pathological conditions, are still in the process of achieving full maturity within the scientific community. This work introduces \lifexep{}, a publicly available software for numerical simulations of the electrophysiology activity of the cardiac muscle, under both physiological and pathological conditions. \lifexep{} employs the monodomain equation to model the heart's electrical activity. It incorporates both phenomenological and second-generation ionic models. These models are discretized using the Finite Element method on tetrahedral or hexahedral meshes. Additionally, \lifexep{} integrates the generation of myocardial fibers based on Laplace-Dirichlet Rule-Based Methods, previously released in Africa et al., 2023, within \lifexfiber{}. 
As an alternative, users can also choose to import myofibers from a file.
This paper provides a concise overview of the mathematical models and numerical methods underlying \lifexep{}, along with comprehensive implementation details and instructions for users. \lifexep{} features exceptional parallel speedup, scaling efficiently when using up to thousands of cores, and its implementation has been verified against an established benchmark problem for computational electrophysiology.  We showcase the key features of \lifexep{} through various idealized and realistic simulations conducted in both physiological and pathological scenarios. Furthermore, the software offers a user-friendly and flexible interface, simplifying the setup of simulations using self-documenting parameter files. \lifexep{} provides easy access to cardiac electrophysiology simulations for a wide user community. It offers a computational tool that integrates models and accurate methods for simulating cardiac electrophysiology within a high-performance framework, while maintaining a user-friendly interface.
\lifexep{} represents a valuable tool for conducting in silico patient-specific simulations.
\end{abstract}


\section{Introduction}\label{sec:background}
Cardiac electrophysiology focuses on the heart conduction system from both the physiological and pathological perspectives: e.g., it involves the study, diagnosis and treatment planning of cardiac arrhythmias \cite{tortora2008, katz2010}.

Nowadays, several clinical tools are widely employed to address these rhythm disorders.
The \ac{ECG} provides a 
recording of the electrical activity of the heart.
Various deviations from sinus rhythm, such as \ac{AF}, \ac{VT}, \ac{LBBB}, can be monitored and identified based on the distinctive shape and morphology of the \ac{ECG} on an individualized patient-specific level \cite{klabunde2011}.
The \ac{ECG} may be combined with other 
imaging data, such as \ac{MRI} and \ac{CT}, or 
electroanatomical maps that are directly recorded on the internal and external surfaces of the heart. Thanks to the aforementioned tools, clinicians can successfully reconstruct both the heart anatomy of a patient and multiple electrophysiology properties of the cardiac tissue.
These information lay the foundations for traditional decision-making in cardiology.
Indeed, medications and surgeries, such as the implantation of pacemakers or cardioverter defibrillators, are planned accordingly~\cite{harrington2011}.

In recent years, the advent of computer models and in-silico simulations in cardiology has enabled the integration of novel scientific tools in standard clinical practice~\cite{quarteroni2019, niederer2019}.
Physics-based mathematical models and data-driven methods are combined to generate digital replicas of human hearts containing a detailed electrophysiology description, both at the cellular level and at the organ scale~\cite{corralacero2020}.
Patient-specific 
clinical data are embedded inside numerical simulations of the cardiac function for precision medicine \cite{peirlinck2021precision}.
For instance, personalized electrophysiology simulations are employed to assess risk stratification of arrhythmias, to define optimal catheter-based ablation targets, and to perform \ac{CRT}~\cite{arevalo2016, prakosa2018, strocchi2020, campos2022}.

In this work, we present \lifexep{}, a publicly available software specifically designed for conducting numerical simulations of cardiac electrophysiology, encompassing both physiological and pathological conditions. \lifexep{} is built upon the foundation of \lifex{} \cite{africa2022lifex}, an open-source, high-performance C++ \ac{FE} numerical solver capable of tackling multi-physics, multi-scale, and multi-domain differential problems. Leveraging the \dealii{} ~\cite{dealII93} \ac{FE} core, \lifex{} was conceived as part of the iHEART project (refer to Section "Funding") with the primary aim of providing the scientific community with a cutting-edge \ac{FE} solver for cardiac modeling.

\begin{figure}[t]
    \centering
    \includegraphics[width=1\textwidth]{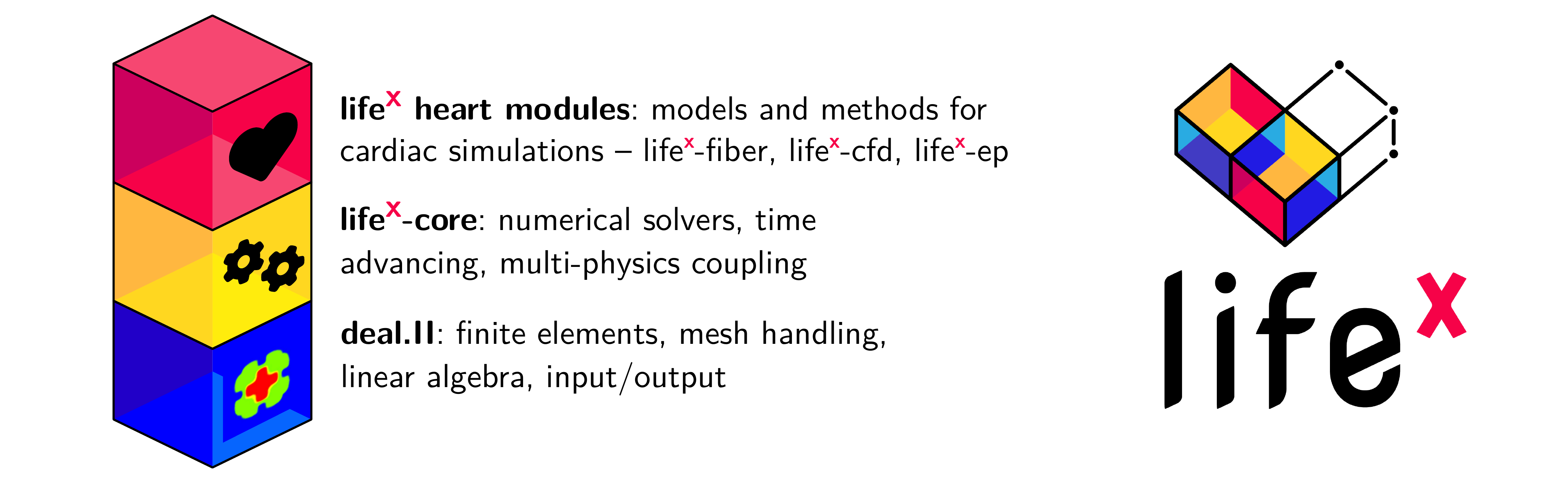}
    \caption{The \lifex{} library encompasses essential features and a framework for numerically solving Finite Element problems. \lifexep{} is a publicly released package designed for cardiac electrophysiology simulations based on \lifex{}. Left picture icons sourced from \url{https://fontawesome.com/license}.}
    \label{fig:lifex_structure}
\end{figure}

As part of the \lifex{} ecosystem, the software released in \lifexep{} has already been widely employed in combination with other modules for cardiac simulations (some of which publicly available: \lifex{} \cite{africa2022lifex}, \lifexfiber{} \cite{lifex-fiber} and \lifexcfd{} \cite{lifex-cfd}, see also Fig.~\ref{fig:lifex_structure}), including electrophysiology \cite{piersanti2021modeling, pagani2021, africa2023matrix, salvador2022role}, mechanics \cite{cicci2022deep, cicci2022efficient, cicci2022projection}, electromechanics \cite{fedele2023comprehensive, salvador2021electromechanical, regazzoni2022cardiac, piersanti20223d}, fluid dynamics \cite{corti2022impact, zingaro2022geometric, zingaro2023electromechanics, fumagalli2022image, marcinno2022computational, bennati2023image, bennati2023turbulence, zingaro2022modeling}, fluid-structure interaction \cite{bucelli2022partitioned}, electro-mechano-fluid interaction \cite{bucelli2022mathematical, bucelli2022stable} and myocardial perfusion \cite{di2022prediction, zingaro2023comprehensive}. This wide range of applications stands as a proof of the flexibility and usability of \lifexep{}. 
Overall, developing a personalized computer model of the human heart electrophysiology that accurately represents all the underlying biological aspects has garnered significant attention. Crucial to achieving this goal is the development of efficient parallel numerical algorithms that can handle computationally intensive tasks with high accuracy. In this context, \lifexep{} proposes a software that offers a user-friendly yet very detailed and highly customizable interface, a variety of modeling and numerical options to choose from, while also ensuring accuracy and computational efficiency, representing a valuable tool for conducting in silico patient-specific simulations.

The present release focuses on the modeling of cardiac electrophysiology. In the following paragraphs, we provide a concise overview of cardiac electrophysiology and the prevailing mathematical methods employed to model it.

\subsection{Cardiac electrophysiology: physiology and modeling}
\label{sec:ep_model}
The heart wall consists of three distinct layers: the internal thin endocardium, the external thin epicardium and the thick muscular cardiac tissue known as the myocardium. The latter is predominantly composed of cardiomyocytes, which are specialized, striated excitable muscle cells responsible for the essential cardiac function.  When these cardiomyocytes are stimulated by an electrical impulse, a change in the electro-chemical balance of the cell membrane results in a series of biochemical reactions that determine a large variation of the \textit{transmembrane potential}, namely the voltage differential between the intra and extracellular spaces of the cell.
The rapid depolarization and subsequent slow repolarization mechanism is known as the \textit{action potential}. This process is triggered and controlled by the opening and closing of voltage-gated ion channels that make the cell membrane permeable to specific ionic species, like sodium, potassium, and calcium. The transmembrane potential changes as a result of the ionic fluxes, which, in turn, are driven by the voltage difference itself.
\begin{figure}[t]
    \centering
    \includegraphics[keepaspectratio, width=0.95\textwidth]{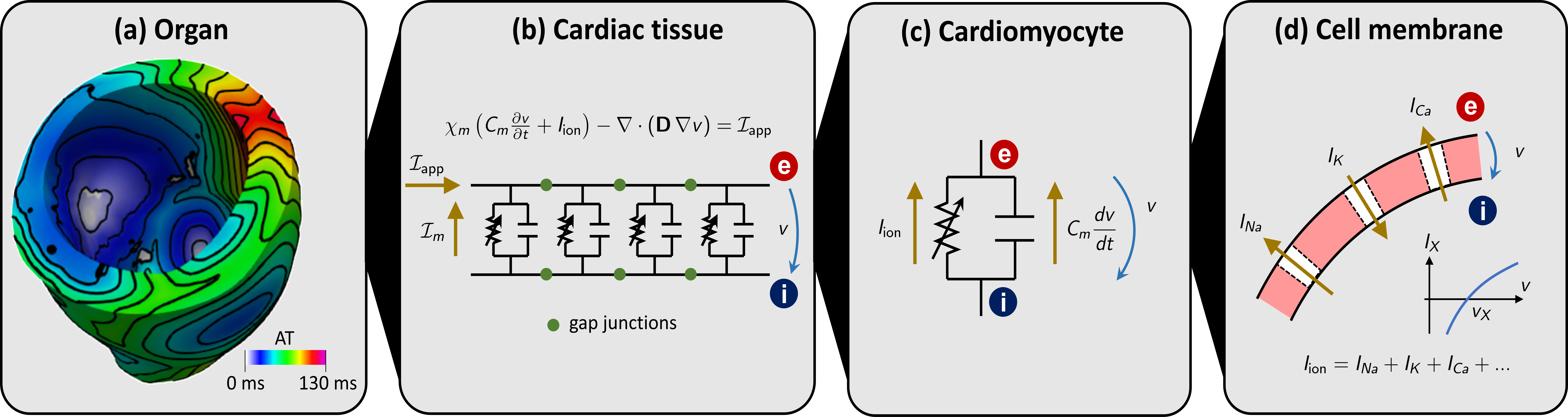}
    \caption{Multiscale cardiac electrophysiology model. From the largest to the smallest spatial
        scale: (a) organ; (b) tissue; (c) cell; (d) membrane.}
    \label{fig:ccs}
\end{figure}
The cells of the heart tissue are connected to each other through \textit{gap junctions}, i.e. intercellular low-resistance ionic channels, which allow the electric signal to travel from cell to cell across the whole cardiac muscle (known as myocardium).
The propagation of the electric signal within the myocardium is highly anisotropic, with preferential directions for conduction determined by the presence of muscle fibers \cite{piersanti2021modeling}.
These fibers are organized in sheets, that determine a second preferential direction of propagation.
For further details on the biophysical mechanisms of cardiac electrophysiology, we refer the interested reader to \cite{tortora2008,katz2010,collifranzone2014}.
A visual representation of the different spatial scales involved in this complex chain of physical processes is given in Fig.~\ref{fig:ccs}.

The mathematical description of these processes is made of two building blocks: a \textit{ionic model}, describing the chemical processes taking place at the cellular scale, and an \textit{action potential propagation model}, representing the spatial propagation of the action potential wavefront at the tissue scale.

\subsubsection{Ionic models}
Several ionic models with different levels of biophysical detail and suitable for specific types of cells have been proposed in the literature \cite{collifranzone2014}. 
Most of them are written in the form of the following \ac{ODEs}
\begin{equation}\label{eqn: ionic model}
    \left\{
    \begin{aligned}
        &\frac{d\Pot(t)}{dt}+\Iion{}(\Pot(t),\Gating{}(t)) = \Iapp(t),  
        & \qquad t \in (0, T],\\
        &\frac{d\Gating{}(t)}{dt}=\RhsGating{}(\Pot(t),\Gating{}(t)), 
        & \qquad t \in (0, T],\\
        &\Pot(0)=\PotStart{}, 
        \, \Gating{}(0)=\GatingStart{}, & 
    \end{aligned}
    \right.
\end{equation}
where the unknowns are $\Pot=\Pot(t)$, the transmembrane potential, and the vector $\Gating{}=\Gating{}(t)=(w_1,\ldots,w_{\NumIonicVariables{}})$, collecting $\NumIonicVariables{}$ ionic variables. 
Ionic variables may include concentrations of different ionic species and gating variables (describing the opening probability of ionic channels).
The term $\Iion{}$ describes the electric current generated by the flux of ionic species across the cell membrane, while $\Iapp$ represents an externally applied current.
The dynamics of ionic variables is modeled 
by the function~$\RhsGating{}$. Different ionic models are characterized by the number of ionic variables $\NumIonicVariables{}$ and the definition of the two functions $\Iion{}$ and $\RhsGating{}$.

The \lifexep{} release includes two phenomenological ionic models (\ac{APf} \cite{aliev1996} and \ac{BO} \cite{buenoorovio2008}) and two physiological ones, one for ventricular cells (\ac{TTP06} \cite{tentusscher2006}) and one for atrial cells (\ac{CRN} \cite{courtemanche1998}).

\subsubsection{Action potential propagation models}

One of the most popular action potential propagation models is the \textit{monodomain equation}, a \ac{PDE} in which the 
spatial propagation of the electric signal through gap junctions is accounted for by a diffusion term \cite{collifranzone2014}.

We consider a computational domain $\Domain{} \subset \mathbf{R}^3$, representing the region of the myocardium of interest (e.g., a slab of cardiac tissue, atria, ventricles or whole heart geometries).
In the whole domain, we consider a local orthonormal triplet of vectors $\fZero$, $\sZero$ and $\nZero$, defining the fiber, the sheetlet and the sheet-normal directions, respectively \cite{piersanti2021modeling}.

A notable feature of \lifexep{} is its capability to divide the domain $\Domain{}$ into a generic number $\NumDomains$ of subdomains, thereby allowing distinct electrophysiology properties to be assigned to each individual region (e.g. different ionic models or different electrical conductivities). In practical implementation, this goal is realized by providing the simulation with a computational mesh, wherein each subdomain is distinctly labeled. More formally, we introduce a partition of $\Domain{}$ into $\NumDomains$ disjoint subdomains, namely $\Domain{1}, \Domain{2}, \dots, \Domain{\NumDomains}$ (more precisely, we assume $\ClosedDomain{} = \cup_{i=1}^{\NumDomains} \ClosedDomain{i}$ and $\Domain{i} \cap \Domain{j} = \emptyset$ for $i \neq j$).
We denote by $\Boundary{i} = \partial\Domain{} \cap \partial\Omega{i}$ the external boundary of the $i$-th subdomain, and by $\Boundary{ij} = \partial\Domain{i} \cap \partial\Omega{j}$ the interface between the $i$-th and the $j$-th subdomains.

For each subdomain $\Domain{i}$, we consider a ionic model, featuring $\NumIonicVariables{i}$ ionic variables, and characterized by the functions $\RhsGating{i}$ and $\Iion{i}$ and initial conditions $\PotStart{i}$ and $\GatingStart{i}$. Moreover, we introduce the unknown $\Gating{i} \colon \Domain{i}\times[0,T] \to \mathbf{R}^{\NumIonicVariables{i}}$.
The transmembrane potential, namely $\Pot \colon \Domain{}\times[0,T] \to \mathbf{R}$, is instead defined in the whole computational domain.

The monodomain model reads as follows \cite{quarteroni2019,collifranzone2014}
\begin{equation}\label{eqn: monodomain}
    \left\{
    \begin{aligned}
        &\frac{\partial \Pot}{\partial t}+\Iion{}(\Pot,\Gating{i})
        -\nabla\cdot(\DiffTens{i}\nabla \Pot)=\Iapp({\bf x},t),
        && \, \mbox{in }\Domain{i}\times(0,T], 
        \\
        &\frac{d\Gating{i}}{dt}=\RhsGating{i}(\Pot,\Gating{i}),
        && \, \mbox{in }\Domain{i}\times(0,T], 
        \\
        & (\DiffTens{i}\nabla \Pot) \cdot \Normal{i} = 0,
        && \, \mbox{on }\Boundary{i}\times(0,T], 
        \\
        & (\DiffTens{i}\nabla \Pot) \cdot \Normal{i} + (\DiffTens{j}\nabla \Pot) \cdot \Normal{j} = 0,
        && \, \mbox{on }\Boundary{ij}\times(0,T], \,
        \\
        & \Pot=\PotStart{i},
        && \, \mbox{in }\Domain{i}\times\{0\}, 
        \\
        & \Gating{i}=\GatingStart{i},
        && \, \mbox{in }\Domain{i}\times\{0\}, 
    \end{aligned}
    \right.
\end{equation}
where $i,j = 1,\dots,\NumDomains$, and the diffusion tensor $\DiffTens{i}$ is defined as
\begin{equation}
    \DiffTens{i} = \sigmal{i} \fZero \otimes \fZero + \sigmat{i} \sZero \otimes \sZero + \sigman{i} \nZero \otimes \nZero,
    \label{eqn: diffusiontensor}
\end{equation}
with  $\sigmal{i}, \sigmat{i}, \sigman{i} \in \mathbf{R}^+$ denoting the longitudinal, transversal and normal conductivities, respectively \cite{piersanti2021modeling}.
Homogeneous Neumann boundary conditions are prescribed on the whole boundary $\partial\Omega$ to impose the condition of electrically isolated domain, with $\Normal{i}$ denoting the outward normal unit vector to the boundary.
On the interface between subdomains, we impose continuity of flux conditions.

The action potential is triggered by an external applied current $\Iapp({\bf x},t)$, that mimics the presence of a natural or artificial pacemaker, or the role of specialized conduction systems, such as the Purkinje network \cite{tortora2008,katz2010} or other relevant bundles.
Since these anatomical entities are not explicitly modeled in this release, \lifexep{} provides great flexibility in the definition of $\Iapp({\bf x},t)$, allowing for the definition of multiple stimuli of various shape and with customized duration and delay.

The possibility of prescribing different properties to different regions of the domain enables several use cases. 
For instance, it allows multi-chamber simulations, in which different ionic models are assigned to atria and ventricles \cite{fedele2023comprehensive}.
Moreover, it also gives the possibility to simulate pathological scenarios, such as the presence of scars or fibrosis \cite{salvador2021electromechanical, frontera2020, frontera2022, pagani2021}. Non-conductive zones can be introduced in the model, by labelling a given subregion of the domain as non-conductive.
In this case, the corresponding portion of the domain will be excluded from the simulation, thus creating a block in the conduction \cite{fedele2023comprehensive}.
Furthermore, grey zones can be modeled as well, by suitably adjusting the tissue conductivities and the ionic model properties in grey zone and fibrotic subregions \cite{pagani2021, salvador2021electromechanical, zahid2016patient, arevalo2016arrhythmia}.

We remark that the monodomain equation is often expressed as
\begin{equation*}
    \Chim\left(\Cm\frac{\partial \Pot}{\partial t}+\IionDIM{}(\Pot,\Gating{})\right) -\nabla\cdot(\DiffTensDIM{}\nabla \Pot)=\Chim\IappDIM({\bf x},t),
\end{equation*}
where $\Cm$ denotes the membrane cell capacitance and $\Chim$ is the membrane surface-to-volume ratio.
This formulation can clearly be recast to the one of \eqref{eqn: monodomain} by the rescaling $\Iion{} = \IionDIM{} / \Cm$, $\Iapp{} = \IappDIM{} / \Cm$ and $\DiffTens{} = \DiffTensDIM{} / (\Chim\Cm)$.
These equivalences are helpful when comparing different literature sources.

\subsection{Cardiac electrophysiology: numerical discretization}
\label{sec:ep_methods}
We partition the temporal domain $(0, T]$ into $N_T$ subintervals with a time step $\Delta t=t_{n+1} - t_{n}$, where $t_n = n\Delta t$ for $n = 0, 1, \dots, N_T$. We denote with a subscript~$n$ the approximation of a variable at time $t_n$ (e.g., $\Pot_n \approx \Pot(t_n)$). Time derivatives are approximated by means of a \ac{BDF} scheme of order $\sigma$, with $\sigma \in \{1, 2, 3\}$, that is
\begin{equation*}
    \frac{\partial f}{\partial t} \approx \frac{\alphabdf f_{n+1} - f\bdf{n}}{\Delta t}\;,
\end{equation*}
where $f\bdf{n}$ is a linear combination of $f_n$, $f_{n-1}$, $\dots$, and $\alphabdf$ is a coefficient depending on the order of the scheme \cite{quarteroni2010numerical}. We also define explicit extrapolations of order $\sigma$, denoted by the subscript EXT (e.g., $f\ext{n+1} \approx f(t_{n+1})$).

The numerical discretization of the system \eqref{eqn: monodomain} is obtained by means of the implicit-explicit time advancing scheme described in \cite{regazzoni2022cardiac}. Given the solution up to time step~$t_n$, to compute the solution at time $t_{n+1}$:
\begin{enumerate}
    \item solve the time-discrete ionic model equations: for $i = 1, \dots, \NumDomains$,
    \begin{equation}\label{eqn: ionic-time-discrete}
        \frac{\alphabdf\Gating{i}_{n+1} - \Gating{i}\bdf{n}}{\Delta t} = \RhsGating{i}(\Pot\ext{n+1},\Gating{i}\ext{n+1},\Gating{i}_{n+1}), 
        \quad \mbox{in }\Domain{i}\;.
    \end{equation}
    We treat some of the ionic variables appearing in $\RhsGating{i}$ with an implicit formulation, and others with an explicit formulation, so that the resulting problem can be solved by directly inverting the equations \cite{regazzoni2022cardiac}.
    \item solve the time-discrete monodomain equation:  for $i,j = 1, \dots, \NumDomains$,
    \begin{equation}\label{eqn: monodomain-time-discrete}
        \left\{\begin{aligned}
            & \begin{aligned}
                \frac{\alphabdf\Pot_{n+1} - \Pot\bdf{n}}{\Delta t} &+\Iion{}(\Pot\ext{n+1},\Gating{i}_{n+1}) \\
                &-\nabla\cdot(\DiffTens{i}\nabla \Pot_{n+1})=\Iapp({\bf x},t),
            \end{aligned}
            && \, \mbox{in }\Domain{i},
            \\
            & (\DiffTens{i}\nabla \Pot_{n+1}) \cdot \Normal{i} = 0,
            && \, \mbox{on }\Boundary{i},
            \\
            & (\DiffTens{i}\nabla \Pot_{n+1}) \cdot \Normal{i} + (\DiffTens{j}\nabla \Pot_{n+1}) \cdot \Normal{j} = 0,
            && \, \mbox{on }\Boundary{ij}.
        \end{aligned}\right.
    \end{equation}
    Notice that the ionic current term is treated explicitly by using the extrapolated potential $\Pot\ext{n+1}$, so that the resulting problem is linear in $\Pot_{n+1}$.
\end{enumerate}

For the spatial discretization, we introduce a tetrahedral or hexahedral mesh over~$\Omega$, and use the \ac{FE} method \cite{quarteroni2017numerical} to approximate the solution variables $\Pot$ and $\Gating{i}$ as piecewise polynomials of order \(p\). \lifexep{} supports polynomials of orders $1$ and $2$ on tetrahedral meshes, and polynomials of arbitrary degree on hexahedral meshes.

The ionic model \eqref{eqn: ionic-time-discrete} is solved at each support point of the degrees of freedom of the \ac{FE} space. The ionic current $\Iion{}$ is also evaluated at every support point, and then interpolated onto quadrature nodes on the interior of the mesh elements, in the approach known as \ac{ICI} \cite{pathmanathan2011significant,rossi2014thermodynamically}. The time-discrete monodomain equation \eqref{eqn: monodomain-time-discrete} is discretized in space using the \ac{FE} method, leading to the following algebraic linear system of equations
\begin{equation}
    \left(\frac{\alphabdf}{\Delta t}\Mass + \Stiff\right) \APot{}_{n+1} = \frac{1}{\Delta t}\Mass{} \APot{}\ext{n+1} - \AIion_{n+1} + \AIapp_{n+1}\;,
    \label{eqn: linear-system}
\end{equation}
where, denoting by $\varphi_j$ the basis functions of the \ac{FE} space, $\Mass$ is the mass matrix of entries $\Mass_{jk} = \int_\Omega \varphi_k \varphi_j d\mathbf x$, $\Stiff$ is the stiffness matrix of entries $\Stiff_{jk} = \sum_{i = 1}^{\NumDomains} \int_{\Domain{i}} \DiffTens{i}\nabla\varphi_k \cdot \nabla\varphi_j d\mathbf x$ and $\AIapp_{n+1}$ is the vector arising from the applied current term $\Iapp$. The vector $\AIion_{n+1}$ arises from the ionic current terms $\Iion{i}$. Exploiting the \ac{ICI} formulation, it can be computed as
\begin{equation*}
    \AIion_{n+1} = \sum_{i = 1}^{\NumDomains} \Mass^i \mathsf{I}_{\text{ion},n+1}^i\;,
\end{equation*}
where $\Mass^i$ is the mass matrix for subdomain $\Domain{i}$, whose entries are $\Mass_{jk}^i = \int_{\Domain{i}} \varphi_k \varphi_j d\mathbf x$, and $\mathsf{I}_{\text{ion},n+1}^i$ is the vector of the evaluations of $\Iion{i}$ at the support points of the \ac{FE} space. The integrals that arise from the \ac{FE} discretization are numerically approximated using the Gauss--Legendre quadrature rule, with the minimum number of points to ensure exact integration of the mass matrix.

\subsection{Overview of existing software}
One of the first landmarks in computational cardiology has been defined by the pioneering work of Hodgkin and Huxley in the mid-1950s \cite{hodgkin1952}.
Since then, a plethora of electrophysiology models has been proposed in the literature \cite{niederer2019,collifranzone2014,lines2003,clayton2020}. These mathematical models feature a different degree of biophysical complexity and act either at the microscopic level, by describing the behavior of single cardiomyocytes, or at the organ scale, where an ensemble of many myocardial cells is considered.

The development of software to perform electrophysiology simulations is still mainly steered by academia and public institutions, despite this variety of mathematical models, the presence of increasingly elaborated numerical methods, and the evolution of computer hardware.
Indeed, several research tools for multi-physics and multi-scale cardiac simulations have been proposed in the last two decades. Among them, an important role is played by \texttt{openCARP} \cite{openCARP}, an open-source C++ simulation environment integrated with \texttt{cellML}, a public repository that encompasses many cell-based mathematical models \cite{cellML}, \texttt{Chaste} \cite{chaste}, an open-source C++ library mainly developed at the University of Oxford, \texttt{Cardioid} \cite{cardioid}, a highly efficient and scalable tool for high resolution electrophysiology simulations mainly developed at the Lawrence Livermore National Laboratory, and \texttt{Alya} \cite{alya}, a high-performance code from the Barcelona Supercomputing Center.
Prominent examples of industrial projects that demonstrate translational research efforts for the heart are given by the Living Heart Project by Dassault Systèmes \cite{livingheartproject2014, livingheartproject2022} and the services and software provided by NumeriCor GmbH\footnote{\url{https://www.numericor.at/}} \cite{vigmond2008solvers} .

The aforementioned tools allow to perform single-chamber, bi-atrial, bi-ventricular or four-chamber heart electrophysiology simulations by means of accurate, yet computationally expensive physiologically-based models, such as the bidomain or monodomain equation coupled with the \ac{TTP06} and \ac{CRN} ionic models \cite{collifranzone2014,tentusscher2006,courtemanche1999}, or the more efficient reaction-eikonal equation~\cite{neic2017}.
Model parameters can be calibrated on a patient-specific basis to match \ac{ECG} or \ac{BSPM}~\cite{gillette2021a}. Different pathological scenarios involving \ac{AF} \cite{loewe2019,azzolin2022}, \ac{VT} \cite{arevalo2016,campos2022} and \ac{LBBB}~ \cite{strocchi2020,gillette2021b} have been addressed by employing these software tools.
Moreover, heterogeneity in the tissue and cellular properties can be prescribed to incorporate the presence of scar, grey zones and fibrosis importing measurements from clinical data, such as \ac{LGE-MRI} \cite{prakosa2018}, contrast enhanced \ac{CT} \cite{sung2020}, or the \acf{IIR} \cite{roney2022}.
This allows to locally vary \ac{CV} and \ac{AP} morphology.
Another application is related to the in-silico assessment of drugs efficacy by means of numerical simulations \cite{margara2021,peirlinck2021}, where model parameters can be tuned to replicate the effects of pharmacological therapies.
Electrophysiology simulations are also used to evaluate gender differences in healthy and pathological conditions involving arrhythmias \cite{peirlinck2021,gonzalezmartin2022}.

Compared to the aforementioned software and libraries, the \lifexep{} solver stands out with several distinctive features and advantages, mostly inherited from the \lifex{} core structure \cite{africa2022lifex}. It is designed to be user-friendly and easy to use, even for biomedical researchers without extensive experience in numerical methods.
The solver is implemented in C++ using advanced programming paradigms and leverages \acf{MPI} for distributed memory parallelism. Moreover, it supports the possibility to import arbitrary meshes with either hexahedral or tetrahedral elements, and incorporates advanced numerical solvers based on the Trilinos linear algebra backend, thus ensuring precise control over the numerical setting and accuracy. The solver also exhibits ideal scalability up to thousands of cores, as demonstrated in Section~\nameref{subsec:scalability}, allowing to efficiently simulate large-scale scenarios. In addition to the numerical and programming features stemming from its foundation on \lifex{}, \lifexep{} offers two options for prescribing myocardial fibers, which can be either imported from a file or generated online taking advantage of the previous release \lifexfiber{} \cite{lifex-fiber}, based on the \acp{LDRBM} presented in \cite{piersanti2021modeling}. Moreover, it supports spatial heterogeneity in the choice of both models and physical coefficients, easily configurable through a convenient parameter file, without the need to access and modify the source code.

In general, \lifexep{} stands out in its ease of use, performance, and compatibility with common I/O (input/output) file formats, as well as its comprehensive and self-contained infrastructure, achieved by combining sophisticated mathematical models with accurate numerical schemes. To the best of our knowledge, none of the packages mentioned above exhibits similar features altogether.

\section{Implementation}\label{sec:implementation}

In this section, we present the technical specifications of \lifexep{} and provide a comprehensive documentation of its user interface. The aim is to guide users through the entire process, from downloading the software to successfully running a full cardiac electrophysiology simulation.

\lifexep{} offers a numerical solver tailored for cardiac electrophysiology, leveraging the mathematical models and numerical algorithms discussed in the previous section. The linear algebra backend is provided by Trilinos \cite{trilinos-website}, integrated into \dealii{}. This incorporates the implementation of various linear solvers (CG and GMRES) and flexible black-box preconditioners (AMG, additive Schwarz, block Jacobi) supported by \lifex{}. The code is inherently parallel and designed to run efficiently on a diverse range of architectures, spanning from personal laptops to High-Performance Computing (HPC) facilities and cloud platforms. To ensure reliability and performance, \lifexep{} has been thoroughly tested on multiple systems, including a cluster node equipped with 192 cores based on Intel Xeon Gold 6238R (2.20 GHz) at MOX, Dipartimento di Matematica, Politecnico di Milano, as well as the \texttt{GALILEO100} supercomputer available at \texttt{CINECA} (Intel CascadeLake 8260, 2.40GHz, technical specifications available at \url{https://wiki.u-gov.it/confluence/display/SCAIUS/UG3.3%3A+GALILEO100+UserGuide}).

For further insights into the core functionalities of \lifex{}, we recommend referring to \cite{africa2022lifex}. Additionally, in the following sections, we provide a concise guide on how to quickly get started and run simulations using \lifexep{}.

\subsection{Running simulations in \texorpdfstring{\lifexep{}}{lifex-ep}}
\label{sec:quick-start}

The distribution and installation process of \lifexep{} is designed to be user-friendly and platform-independent. The software is conveniently provided as a binary \texttt{AppImage}\footnote{\url{https://appimage.org/}} executable, which can be obtained from \zenodourl{}. Along with the executable, all the necessary input files required to reproduce the numerical results presented in the subsequent sections are included.

The adoption of the \texttt{AppImage} format ensures a universal package that is compatible with \texttt{x86-64 Linux} operating systems, eliminating the need for multiple distribution-specific versions. From the user's perspective, this translates to a seamless and straightforward \textit{download-then-run} experience, without the hassle of manually managing system dependencies.

As an \texttt{AppImage}, \lifexep{} has been built on \texttt{Debian Buster}\footnote{\url{https://www.debian.org/releases/}}, which corresponds to the current \texttt{oldoldstable} version. This follows the principle of \textit{``Build on old systems, run on newer systems"}\footnote{\url{https://docs.appimage.org/introduction/concepts.html}}. Consequently, the software is expected to function on virtually any recent \texttt{x86-64 Linux} distribution, provided that \texttt{glibc}\footnote{\url{https://www.gnu.org/software/libc/}} version \texttt{2.28} or higher is installed.

Once downloaded and extracted the \lifexep{} archive, the \texttt{AppImage} file needs to be made executable by typing the following command in a terminal:
\begin{lstlisting}[language=bash]
$ chmod +x lifex_electrophysiology-1.5.0-x86_64.AppImage
\end{lstlisting}

Finally, \verb|lifex_electrophysiology-1.5.0-x86_64.AppImage| can be executed with:
\begin{lstlisting}[language=bash]
$ ./lifex_electrophysiology-1.5.0-x86_64.AppImage [ARGS]...
\end{lstlisting}

No root permissions are required for the commands mentioned above to run successfully. However, it is important to note that the \texttt{AppImage} relies on the userspace filesystem framework called \texttt{FUSE}\footnote{\url{https://www.kernel.org/doc/html/latest/filesystems/fuse.html}}. Please ensure that \texttt{FUSE} is installed on your system. If you encounter any errors, the following commands may be helpful in resolving the issue:
\begin{lstlisting}[language=bash]
$ ./lifex_electrophysiology-1.5.0-x86_64.AppImage \
      --appimage-extract
    
$ squashfs-root/usr/bin/lifex_ep [ARGS]...
\end{lstlisting}
Additionally, we recommend referring to the \texttt{AppImage} troubleshooting guide\footnote{\url{https://docs.appimage.org/user-guide/troubleshooting/fuse.html}}.

The following command will provide an inline help that includes detailed information about all the available command line options and their purpose:
\begin{lstlisting}[language=bash]
$ ./lifex_electrophysiology-1.5.0-x86_64.AppImage -h
\end{lstlisting}

The executable allows to run test cases with an arbitrary number of disjoint subdomains \(\Omega_i,\,i=1,\dots,N\), which are also referred to as \textit{volumes}. The configuration of the simulation is supplied through a parameter file. The user can generate a template parameter file using the following command:
\begin{lstlisting}[language=bash]
$ ./lifex_electrophysiology-1.5.0-x86_64.AppImage \ 
      -g [minimal|full] \
      [-vol <volume labels>...]
\end{lstlisting}
To match different user needs, the level of detail in the parameter files can be adjusted using the optional \texttt{minimal} or \texttt{full} option after the \texttt{-g} flag. The \texttt{minimal} option reduces the level of detail, making it suitable for initial usage of \lifexep{}. On the other hand, the \texttt{full} option increases the level of detail, exposing advanced options, such as parameter choices
and detailed options on linear algebra and preconditioning, among others. If the user does not specify either the \texttt{minimal} or \texttt{full} option, an intermediate verbosity level is selected by default. In all cases, parameters that are not present in the file will retain their default values. This flexibility allows users to customize the level of detail in the parameter files according to their specific requirements and familiarity with the software.

The parameters are written in a plain text file, organized as a list of key-value pairs grouped in subsections, which describe the configuration for the simulation to be run. Each parameter is accompanied by a brief documentation within the parameter file itself, explaining its meaning.

In the provided command, the optional argument \texttt{-vol <volume labels>} allows to specify a list of user-defined labels. These labels are used to differentiate each subdomain, enabling the selection of heterogeneous model options such as ionic model type, coefficients, and electrical conductivities for each subdomain. If not specified, a single subdomain characterized by its global \verb|Volumetric parameters| is assumed to exist.

The following example illustrates a parameter file specifying three subdomains: \texttt{Healthy}, \texttt{Fibrosis}, and \texttt{Scar}:

\begin{lstlisting}[language=prm]
subsection Electrophysiology
  # ...
  subsection Physical constants and models
    # ...
    subsection Healthy
      set Material IDs   = 1
      set Ionic model    = TTP06
      # ...
      subsection Ionic model parameters
        # ...
      end
    end
    
    subsection Fibrosis
      set Material IDs   = 2 3
      set Ionic model    = Bueno-Orovio
      # ...
      subsection Ionic model parameters
      # ...
      end
    end
    
    subsection Scar
      set Material IDs   = 4
    # ...
    end
  # ...
  end
# ...
end
\end{lstlisting}

The input mesh is expected to include at least one volumetric tag corresponding to each of the subdomains that need to be differentiated. These subdomain tags are specified in the \verb|Material IDs| list, located under each respective subdomain section.

The parameter file includes a dedicated section called \verb|Fiber generation|, which serves the purpose of enabling the importing from a file of the myocardial fibers or the online generation of them on various geometry types, such as slabs, ventricles, and atria. To accomplish this, \lifexep{} incorporates the functionalities of its predecessor, \lifexfiber{} \cite{lifex-fiber}, which utilizes the \aclp{LDRBM} presented in \cite{piersanti2021modeling}. This unique feature of \lifexep{} sets it apart from other existing software alternatives.

Once the user has edited the parameter file, the simulation can be started using the following command:
\begin{lstlisting}[language=bash]
$ ./lifex_electrophysiology-1.5.0-x86_64.AppImage \ 
      [-vol <volume labels>...] \
      [-f parameter_file_name.prm] \
      [-o output_folder]
\end{lstlisting}
When executing the command, the volume labels provided must match the ones used for generating the parameter file. It is essential to use consistent volume labels throughout the process to ensure proper identification and configuration of the subdomains within the simulation.

A parallel simulation is started prepending the command with the \texttt{mpirun} or \texttt{mpiexec} wrapper (which may vary depending on the MPI implementation available), e.g.:
\begin{lstlisting}[language=bash]
$ mpirun -n N ./lifex_electrophysiology-1.5.0-x86_64.AppImage \ 
  ...
\end{lstlisting}
where \texttt{N} represents the desired number of parallel processes. The binary package supports parallel execution using \texttt{MPICH} (\url{https://www.mpich.org/}) version \texttt{4.0} or higher.

The parameter file also includes options that enable the serialization of the solution, allowing the simulation to be paused or stopped at any point and then resumed at a later time using the serialized data. This feature is particularly useful when dealing with long-running simulations or when unexpected interruptions occur.

\subsection{License and third-party software}
This work is is copyrighted by the \lifexep{} authors and licensed under the Creative Commons Attribution Non-Commercial No-Derivatives 4.0 International License\footnote{\url{http://creativecommons.org/licenses/by-nc-nd/4.0/}}.

It should be noted that \lifexep{} incorporates several third-party libraries, which are separately copyrighted by their respective authors and whose use is covered by various permissive licenses.

Third-party software bundled with (in binary form),
required by, copied, modified, or explicitly used in
\lifexep{} include:
\begin{description}
    \item[\lifex{}\footnote{\url{https://lifex.gitlab.io/}} \cite{africa2022lifex}:] the open-source, high-performance software providing the core functionalities for the numerical solution of the \acl{FE} problems described in the previous section;
    \item[\texttt{deal.II}\footnote{\url{https://www.dealii.org/}} \cite{dealII93}:] it provides support to mesh handling, assembling and solving \acl{FE} problems (compiled with enabled support to \texttt{Trilinos}\footnote{\url{https://trilinos.github.io/}} for linear algebra data structures and solvers) and to input/output functionalities;
    \item[\texttt{Boost}\footnote{\url{https://www.boost.org/}} \cite{schaling2011boost}:] its modules \texttt{Filesystem} and \texttt{Math} are used for manipulating files/directories and for advanced mathematical functions and interpolators, respectively;
    \item[\texttt{VTK}\footnote{\url{https://vtk.org/}} \cite{schroeder2006visualization}:] it is used for importing external surface or volume input data and coefficients appearing in the mathematical formulation.
\end{description}
Some of the packages listed above, as stated by their respective authors, rely on additional third-party dependencies that may also be bundled (in binary form) with
\lifexep{}, although not used directly.
These dependencies include:
\texttt{ADOL-C}\footnote{\url{https://github.com/coin-or/ADOL-C}},
\texttt{ARPACK-NG}\footnote{\url{https://github.com/opencollab/arpack-ng}},
\texttt{BLACS}\footnote{\url{https://www.netlib.org/blacs/}},
\texttt{Eigen}\footnote{\url{https://eigen.tuxfamily.org/}},
\texttt{FFTW}\footnote{\url{https://www.fftw.org/}},
\texttt{GLPK}\footnote{\url{https://www.gnu.org/software/glpk/}},
\texttt{HDF5}\footnote{\url{https://www.hdfgroup.org/solutions/hdf5/}},
\texttt{HYPRE}\footnote{\url{https://www.llnl.gov/casc/hypre/}},
\texttt{METIS}\footnote{\url{http://glaros.dtc.umn.edu/gkhome/metis/metis/overview}},
\texttt{MPICH}\footnote{\url{https://www.mpich.org/}},
\texttt{MUMPS}\footnote{\url{http://mumps.enseeiht.fr/index.php?page=home}},
\texttt{NetCDF}\footnote{\url{https://www.unidata.ucar.edu/software/netcdf/}},
\texttt{OpenBLAS}\footnote{\url{https://www.openblas.net/}},
\texttt{PETSc}\footnote{\url{https://www.mcs.anl.gov/petsc/}},
\texttt{ParMETIS}\footnote{\url{http://glaros.dtc.umn.edu/gkhome/metis/parmetis/overview}},
\texttt{ScaLAPACK}\footnote{\url{https://www.netlib.org/scalapack/}},
\texttt{Scotch}\footnote{\url{https://gitlab.inria.fr/scotch/scotch}},
\texttt{SuiteSparse}\footnote{\url{https://people.engr.tamu.edu/davis/suitesparse.html}},
\texttt{SuperLU}\footnote{\url{https://portal.nersc.gov/project/sparse/superlu/}},
\texttt{oneTBB}\footnote{\url{https://oneapi-src.github.io/oneTBB/}},
\texttt{p4est}\footnote{\url{https://www.p4est.org/}}.

These libraries are free software and have relatively few restrictions on their use. However, please note that different terms may apply. For detailed information on the licenses and copyright statements for these packages, please refer to the content of the folder \verb|share/doc/licenses/| in the \lifexep{} archive.

\section{Results and Discussion}\label{sec:results}
To highlight the versatility of \lifexep{}, this section presents the results obtained in various numerical examples conducted on a range of idealized and realistic geometries, encompassing both physiological and pathological scenarios. To facilitate the reproduction of these test cases, this section provides a \textit{getting started} guide, detailing the following steps:
\begin{itemize}
    \item generating and importing the input data (e.g., computational meshes and fibers);
    \item configuring the parameter files specific to each test case and executing the corresponding simulation;
    \item performing post-processing of the results and visualizing the output.
\end{itemize}
Finally, a strong scalability test is presented in Section~\nameref{subsec:scalability}.

\subsection{Input data}
\label{subsec:data}
All parameter files and meshes related to the numerical simulations described below can be downloaded from the \lifexep{} release archive \zenodourl{}.

\begin{figure}[t]
    \centering
    \includegraphics[keepaspectratio, width=0.95\textwidth]{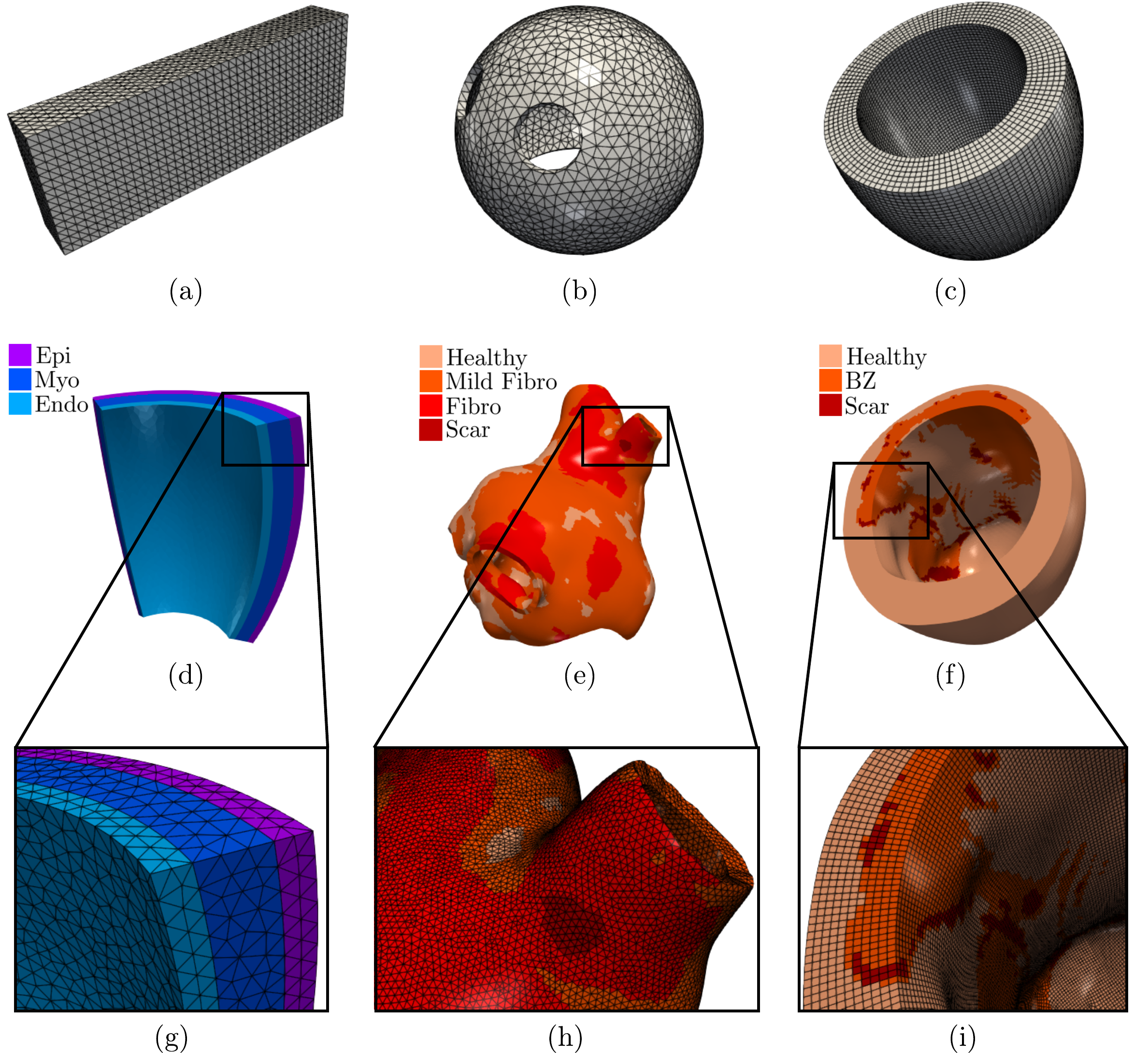}
    \caption{Domains and meshes used in the numerical examples. Domains in the top row (a-c) are composed of a single subdomain, while for domains in the middle row (d-f), zoomed on the bottom (g-i), colors are used to differentiate the subdomains.}
    \label{fig:mesh}
\end{figure}

This guide provides various pre-configured hexahedral and tetrahedral meshes, including:
\begin{itemize}
    \item four idealized geometries: cardiac slab tissue, left atrium, left ventricle, and ventricular slab, see~\cref{fig:mesh}(a-d,~g);
    \item two realistic geometries: left atrium and left ventricle, see~\cref{fig:mesh}(e-f, h-i).
\end{itemize}
We emphasize that the provided example meshes are solely used to illustrate the \lifexep{} features, as users have the flexibility to input any (idealized or patient-specific) meshes into \lifexep{}.

The cardiac tissue slab, idealized left atrial, and left ventricular meshes are characterized by a single volumetric tag representing the entire myocardium, \cref{fig:mesh}(a-c). On the other hand, the ventricular slab, realistic left atrial, and left ventricular meshes have multiple distinct volumetric tags,~\cref{fig:mesh}(d-i). The ventricular slab is divided into sub-endocardial, myocardial, and sub-epicardial layers,~\cref{fig:mesh}(d, g). The realistic left atrium and left ventricle include regions with scars, grey zones, and fibrosis,~\cref{fig:mesh}(e, f, h,  i). The former meshes are used for \textit{single volume simulations}, while the latter for \textit{multi-volume simulations}. In all cases, volumetric tags must be defined during the mesh generation process \cite{fedele2021polygonal}.

The geometrical models for the idealized cardiac meshes were created using the built-in CAD engine of \texttt{gmsh}\footnote{\url{https://gmsh.info/doc/texinfo/gmsh.html}}, an open-source 3D \ac{FE} mesh generator. The \texttt{gmsh} scripts used to generate these meshes are included in the \lifexep{} release archive (\zenodourl{}). For detailed information on the 
syntax of the scripts, we refer to the online documentation of \texttt{gmsh}. Tetrahedral mesh generation for the slab tissue, idealized left atrium, and ventricular slab are also performed using \texttt{gmsh}. On the other hand, the idealized and realistic left ventricular hexahedral meshes were instead generated using the \texttt{cubit}\footnote{\url{https://cubit.sandia.gov/}} mesh generation software toolkit. Finally, the realistic left atrial mesh was 
perfomed using the semi-automatic meshing tools developed in \cite{fedele2021polygonal}, based on the Vascular Modelling Toolkit (\texttt{vmtk}) software \cite{antiga2008vascular}.

The realistic left atrium and left ventricle, also containing the scar and fibrotic regions, have been produced starting from the openly available meshes used in \cite{roney2022predicting} (for the left atrium\footnote{\url{https://zenodo.org/record/5801337}}) and in \cite{mendoncacosta2021} (for the left ventricle\footnote{\url{https://kcl.figshare.com/articles/dataset/A_Virtual_Cohort_of_Twenty-four_Left-ventricular_Models_of_Ischemic_Cardiomyopathy_Patients/16473903}}). For the latter, we used the (\texttt{vmtk}) software \cite{antiga2008vascular}, in conjunction with the \texttt{cubit} mesh generator.
Mesh files can be specified in the \lifexep{} parameter file for electrophysiology simulation within the section named \texttt{Mesh and space discretization}, by setting the proper element type (tetrahedra or hexahedra), the polynomial order to be used in the \ac{FE} discretization and an appropriate mesh rescaling factor (e.g., if the coordinates in the input mesh file are given in millimeters, \texttt{Scaling factor~=~1e-3}, since \lifexep{} internally handles all physical quantities in the International System of Units.

\begin{lstlisting}[language=prm]
subsection Electrophysiology
  subsection Mesh and space discretization
    # Available options are Hex for hexahedra 
    #                   and Tet for tetrahedra
    set Element type    = Hex
    set FE space degree = 1
    subsection File
      set Filename       = /path/to/mesh/mesh.msh
      set Scaling factor = 1e-3 # [mm] to [m]
    end
  end
# ...
end
\end{lstlisting}

Regarding the prescription of the myocardial fiber architecture, an essential building block for cardiac electrophysiology simulations, \lifexep{} provides two options. Users can either import the myofibers from a file or generate them online fusing \acp{LDRBM}~\cite{piersanti2021modeling} by incorporating the \lifexfiber{} release package, which was recently published in \cite{lifex-fiber}.

To generate the myocardial fibers using \acp{LDRBM}, users can select the appropriate \texttt{Geometry type} within the parameter file under the section named \texttt{Fiber generation}. Each \texttt{Geometry type} corresponds to the specific \ac{LDRBM} applicable for different geometries, such as (ventricular and spherical) slabs, (cut at base and complete) left ventricular and left atrial geometries.  The parameters related to the fiber generation for a particular \texttt{Geometry Type} are located within a subsection with the same name. We refer to \cite{lifex-fiber} for a more comprehensive guide on configuring \acp{LDRBM}, and to \cite{piersanti2021modeling} for an in-depth mathematical description.

\begin{lstlisting}[language=prm]
subsection Fiber generation
  subsection Mesh and space discretization
    # Available options are Import from file, Slab, 
    # Left ventricle, Left ventricle complete, 
    # Left atrium
    set Geometry type = Slab
  end
  # ...
  subsection Slab
  # ...
  end
end
\end{lstlisting}
\begin{figure}[t]
    \centering
    \includegraphics[keepaspectratio, width=0.95\textwidth]{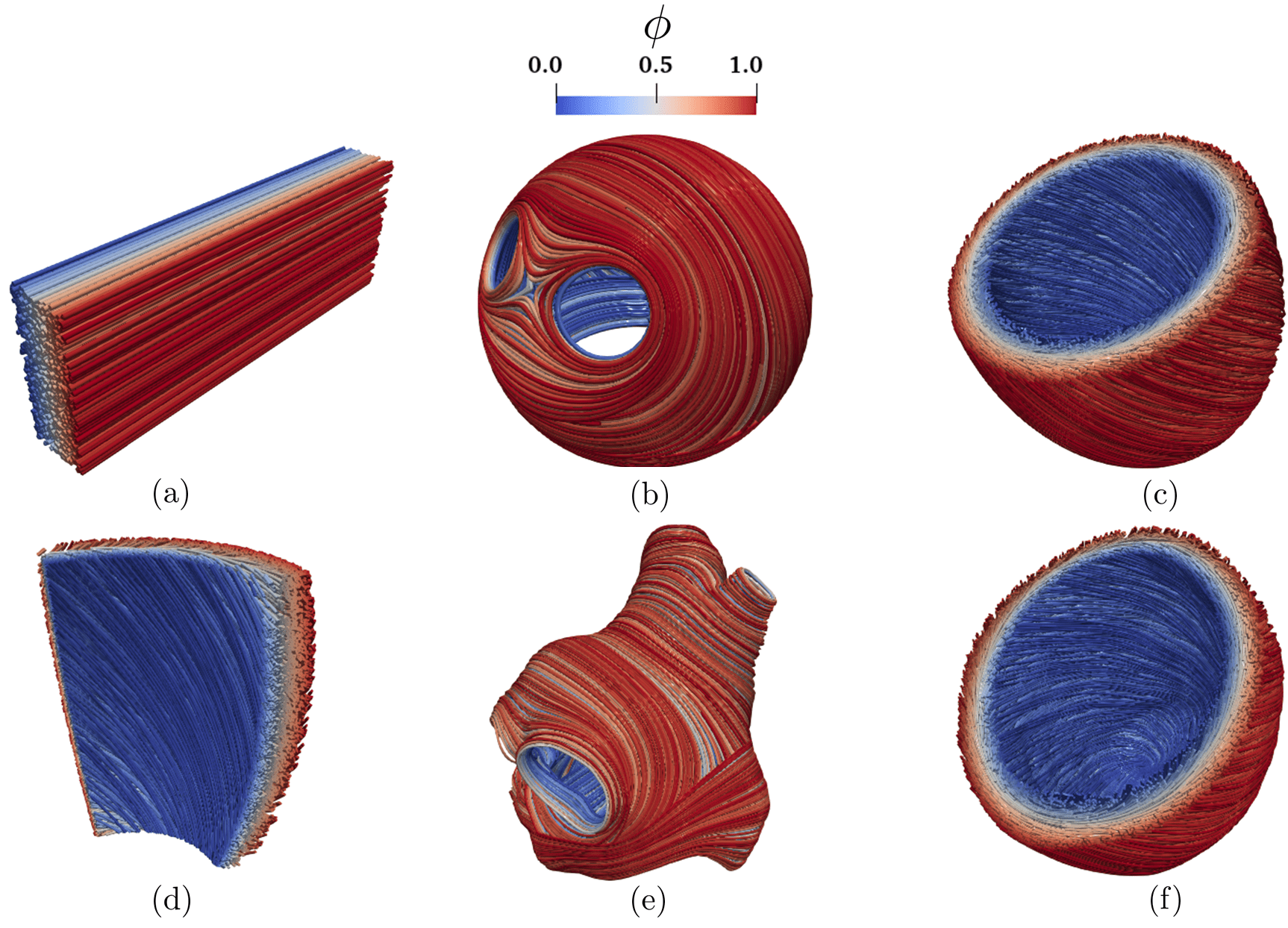}
    \caption{Fiber field computed using \aclp{LDRBM} \cite{lifex-fiber,piersanti2021modeling} visualized as streamlines: (a) Slab tissue; (b) Idealized left atrium; (c) Idealized left ventricle; (d) Ventricular slab; (e) Realistic left atrium; (f) Realistic left ventricle. The Laplace solution $\phi$ is the transmural function where $\phi=0$ on the endocardium and $\phi=1$ on the epicardium.}
    \label{fig:fiber}
\end{figure}
As an alternative, users can choose to import myofibers from a file by setting \texttt{Geometry type = Import from file}. The myofibers are imported from a \texttt{VTU} file format with unstructured grid data. In this file, the three fiber directions ($\fZero$, $\sZero$, and $\nZero$ representing the fiber, sheet, and sheet-normal directions, respectively) must be normalized and embedded as point-data arrays. Additionally, users need to set an appropriate geometry rescaling factor (e.g., if the coordinates in the input fiber file are given in millimeters, \texttt{Scaling factor = 1e-3}). This ensures that the imported myofibers align properly with the geometry of the cardiac electrophysiology model.

\begin{lstlisting}[language=prm]
subsection Fiber generation
  subsection Mesh and space discretization
    set Geometry type = Import from file
  end
  # ...
  subsection Import fibers from file
    # VTU file containing f0, s0, n0.
    set VTU filename = /path/to/myofibers_to_import.vtu
    # Comma-separated list of array names for f0, s0, n0.
    set Array names  = fiber, sheet, sheet_normal
    set Geometry scaling factor = 1e-3 # [mm] to [m].
  end
# ...
end
\end{lstlisting}

The myocardial fibers, generated using \ac{LDRBM}, for the electrophysiology simulations presented hereafter are illustrated in~\cref{fig:fiber}.
\subsection{Physiological electrophysiology simulations}
\label{subsec:pyshio}
We present the physiological electrophysiology simulations applied to the set of idealized geometries, namely a rectangular slab of cardiac tissue, an idealized left atrium, an idealized left ventricle and a layered ventricular slab.

\subsubsection{Slab benchmark}
\label{subsec:slab}
To perform software verification, we consider the N-version slab benchmark proposed in \cite{niederer2011verification}. This benchmark involves a rectangular slab  domain of size $(3 \times 7 \times 20) \times 10^{-3}$ m, as depicted in Figure \ref{fig:mesh}(a). The fiber directions within the domain are oriented along the long axis (0.02 m), as shown in Figure \ref{fig:fiber}(a). For comprehensive modeling and geometrical information regarding the benchmark definition, we refer to the original paper \cite{niederer2011verification}.

The simulation can be run by first generating the parameter file using
\begin{lstlisting}[language=prm]
$ ./lifex_electrophysiology-1.5.0-x86_64.AppImage -g
\end{lstlisting}
then configuring the simulation by editing the \texttt{lifex\_electrophysiology.prm} default parameter file, and finally running the simulation by typing
\begin{lstlisting}[language=prm]
$ ./lifex_electrophysiology-1.5.0-x86_64.AppImage \
      -f lifex_electrophysiology.prm
\end{lstlisting}
The same can be obtained by directly using the already prepared parameter file \texttt{params\_slab.prm} uploaded in the release archive:
\begin{lstlisting}[language=prm]
$ ./lifex_electrophysiology-1.5.0-x86_64.AppImage \ 
      -f params_slab.prm
\end{lstlisting}

To match with \cite{niederer2011verification}, precise specifications are provided for the ionic model employed (in this test case, the \ac{TTP06} model \cite{tentusscher2006}), including initial conditions and conductivity values along the myofiber directions (longitudinal, transversal, and normal).
\begin{lstlisting}[language=prm]
subsection Electrophysiology
  # ...
  subsection Physical constants and models
  # ...
    subsection Volumetric parameters
      set Ionic model = TTP06
      subsection Monodomain conductivities
        set Longitudinal conductivity = 0.95298e-4
        set Transversal conductivity  = 0.12576e-4
        set Normal conductivity       = 0.12576e-4
      end
      # ...
      subsection Ionic model parameters
        # ...
        subsection TTP06
          set Cell type                 = Epicardium
          subsection Physical constants
            # ...
            subsection Initial conditions
              set Transmembrane potential = -85.23e-3
              set M                       = 0.00172
              set H                       = 0.7444
              set J                       = 0.7045
              set Xr1                     = 0.00621
              set Xr2                     = 0.4712
              set Xs                      = 0.00095
              set S                       = 0.999998
              set R                       = 2.42e-8
              set D                       = 3.373e-5
              set F                       = 0.7888
              set F2                      = 0.9755
              set FCass                   = 0.9953
              set Cai                     = 0.000126
              set CaSR                    = 3.64
              set CaSS                    = 0.00036
              set Nai                     = 8.604
              set Ki                      = 136.89
              set RR                      = 0.9073
            end
          end
        end
      end
    end
  end
  # ...
end
\end{lstlisting}
The applied stimulus current, delivered to a volume of $(1.5 \times 1.5 \times 1.5) \times 10^{-3}$ m, situated at one corner of the slab, is prescribed in the subsection \texttt{Applied current}. The stimulus has a duration of $2 \times 10^{-3}$ s and an amplitude of $35.714$ V/s.
\begin{lstlisting}[language=prm]
subsection Electrophysiology
  # ...
  subsection Applied current
    # ...
    subsection Cubic
      set Active                = true
      set Impulse sites         = 0.00075 0.00075 0.00075
      set Impulse amplitudes    = 35.714
      set Impulse length        = 1.5e-3
      set Impulse initial times = 0e-3
      set Impulse durations     = 2e-3
    end
  end
  # ...
end
\end{lstlisting}
\begin{figure}[t!]
	\centering
	\includegraphics[keepaspectratio, width=0.95\textwidth]{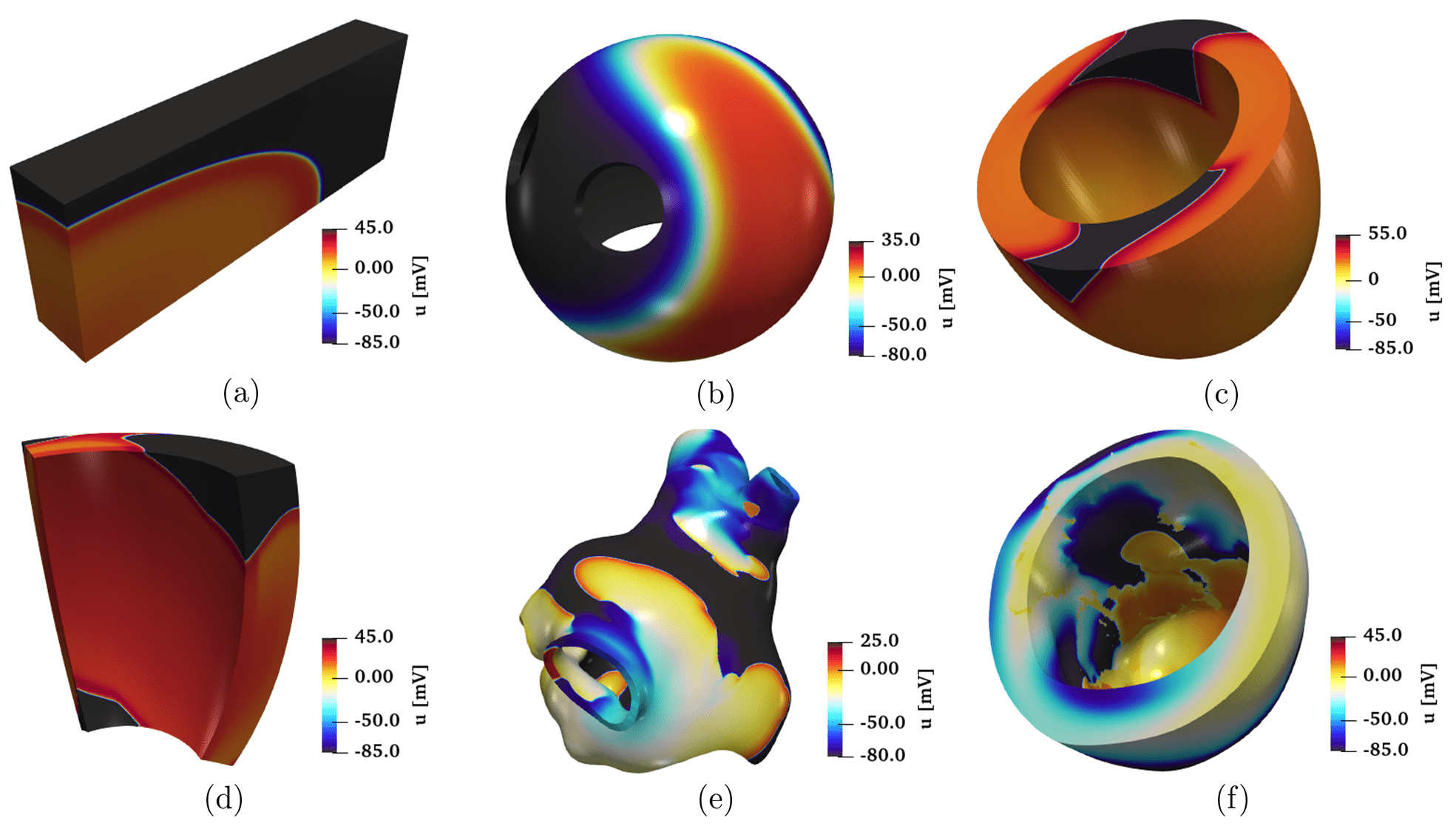}
	\caption{Snapshots of the transmembrane potential for all test cases: (a) Slab tissue; (b) Idealized left atrium; (c) Idealized left ventricle; (d) Ventricular slab; (e) Realistic left atrium; (f) Realistic left ventricle.}
	\label{fig:u}
\end{figure}
\begin{figure}[ht!]
	\centering
	\includegraphics[keepaspectratio, width=0.95\textwidth]{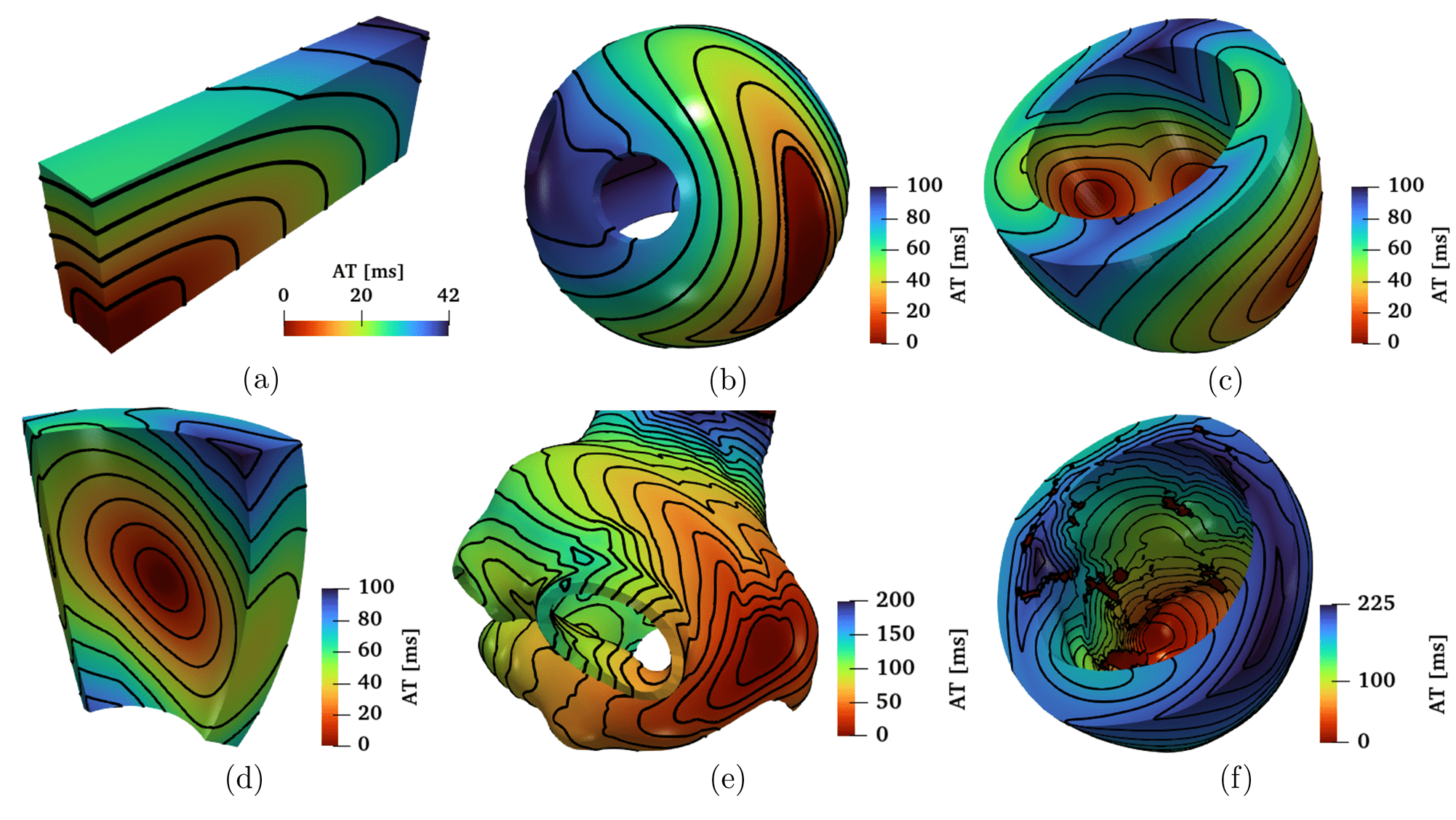}
	\caption{Activation maps computed for all the test cases: (a) Slab tissue; (b) Idealized left atrium; (c) Idealized left ventricle; (d) Ventricular slab; (e) Realistic left atrium; (f) Realistic left ventricle.}
	\label{fig:at}
\end{figure}
In \cref{fig:u}(a) and \cref{fig:at}(a), we display a snapshot of the transmembrane potential and the total activation time computed as output of the
numerical simulation, respectively. The computation of the activation time is evaluated in a given point in the cardiac muscle as the time when the transmembrane potential derivative~$\frac{\partial u}{\partial t}$ reaches its maximum value. This can be enabled in the parameter file under the subsection named \texttt{Activation time}.

\begin{figure}[t!]
    \centering
    \includegraphics[keepaspectratio, width=0.95\textwidth]{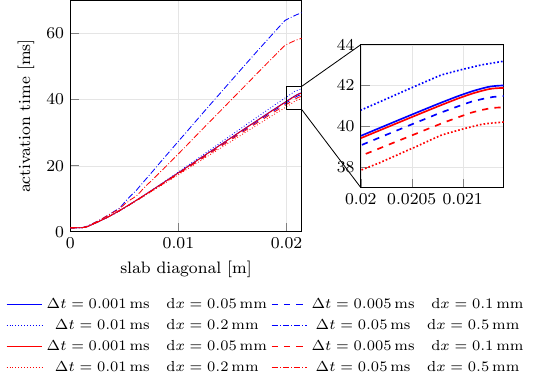}
    \caption{Activation times evaluated along the cuboid diagonal line in the N-version benchmark problem \cite{niederer2011verification} for all the numerical solutions performed with \lifexep{} at different refinements in space and time. Red lines=Hexahedral simulations; Blue lines=Tetrahedral simulations.}
    \label{fig:convergence}
\end{figure}

\begin{figure}[h!]
    \centering
    \includegraphics[keepaspectratio, width=0.95\textwidth]{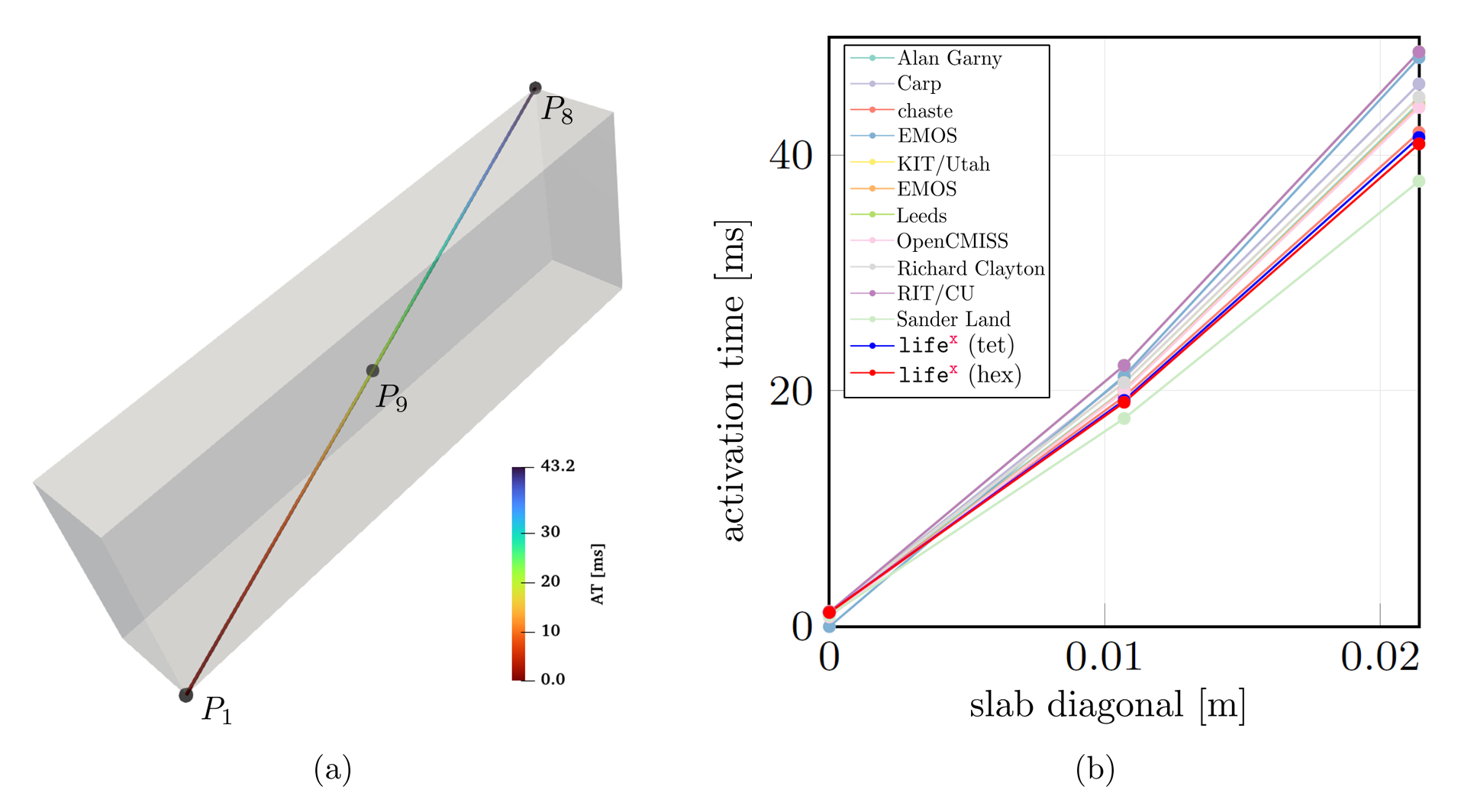}
    \caption{(a) Activation times evaluted along the cuboid diagonal and in the points P1, P9 and P8 in the N-version benchmark problem \cite{niederer2011verification}. (b) Comparison of the \lifexep{} numerical solutions (with Red line=Hexahedral mesh and Blue line=Tetrahedral mesh) with respect to the other codes partecipating to the benchmark problem \cite{niederer2011verification}.}
    \label{fig:verification}
\end{figure}

The problem was solved using both tetrahedral and hexahedral conforming meshes and the \ac{BDF}2 scheme, using eight combinations of spatial resolutions ($dx = [0.5, 0.2, 0.1, 0.05] \times 10^{-3}$ m) and time steps ($\Delta t = [0.05, 0.01, 0.005, 0.001] \times 10^{-3}$~s), as reported in~\cref{fig:convergence}. The results showcase the activation time at points along the diagonal line of the slab, as also shown in Figure \ref{fig:verification}(a). For both the tetrahedral and hexahedral meshes, the numerical solutions converge towards the finer space-time discretization ($\Delta t = 0.001 \times 10^{-3}$ s, $dx = 0.05 \times 10^{-3}$ m), yielding a latest activation time of $41.8$ ms and $42.0$ ms, respectively. These findings align with the values reported in the original N-version benchmark paper \cite{niederer2011verification}.

To further evaluate the \lifexep{} results in comparison to the other eleven codes participating in the N-version benchmark \cite{niederer2011verification}, we report in \cref{fig:verification} the activation time (for $\Delta t = 0.005 \times 10^{-3}$ s, $dx = 0.1 \times 10^{-3}$ m) at specific points along the slab diagonal (namely P1-P8-P9, see~\cref{fig:verification}(a)) for all the codes, including \lifexep{}. The activation times at points P1-P8-P9 for all the other eleven codes are available in the electronic supplementary material of~\cite{niederer2011verification}. As shown in \cref{fig:verification}(b), the \lifexep{} results fall within the range of activation time values computed by the other codes. Moreover, at such refinement level ($dt = 0.005 \times 10^{-3}$ s, $dx = 0.1 \times 10^{-3}$ m), the \lifexep{} numerical solutions (in both tetrahedral and hexahedral meshes) are in excellent agreement with the reference converged activation time values of $42$-$43$ ms computed at point P8, where the accumulation of errors tends to occur as the wave propagates across the cuboid.

\subsubsection{Idealized left atrium}
\label{subsec:ideal_la}
We simulate the propagation of the electrophysiology wavefront in an idealized left atrial geometry (see \cref{fig:mesh}(b)) using the \ac{APf} ionic model \cite{aliev1996}. To run the simulation, the user can modify the generic parameter file 
\newline
\texttt{lifex\_electrophysiology.prm} or use the predefined file
\texttt{params\_ideal\_la.prm} available in the release archive.

\begin{lstlisting}[language=prm]
    $ ./lifex_electrophysiology-1.5.0-x86_64.AppImage \
          -f params_ideal_la.prm
\end{lstlisting}
We employ a second-order \ac{BDF} temporal scheme with a time step $\Delta t=10^{-4}$ s and a final time $T=0.15$ s. Furthermore, we initialize the ionic model by running a 1000-cycle long single-cell simulation applying a stimulus with period $0.8$ s.
\begin{lstlisting}[language=prm]
subsection Electrophysiology
  # ...
  subsection Physical constants and models
    set Time 0D simulation for variables initialization = 800
    subsection Volumetric parameters
      set Ionic model = Aliev-Panfilov
      # ...
      subsection Ionic model parameters
        # ...
        subsection Time solver 0D
          set Time step = 1e-4
        end
        subsection Applied current 0D
          set Initial times = 0.0
          set Durations     = 4e-3
          set Amplitudes    = 1.1628e3
          set Period        = 0.8
        end
        # ...
      end
    end
  end
  # ...
end
\end{lstlisting}
The simulation is initiated using a single spherical impulse with a radius of $r=3~\times~10^{-3}$ m. The parameters related to the impulse site, amplitude, duration, and initial time are specified in the \texttt{Applied current/Spherical} subsection of the parameter file. Notice that, when utilizing the \ac{APf} model, it is essential to properly rescale the impulse amplitude and also conductivity values embedded within the diffusion tensor \eqref{eqn: diffusiontensor} in accordance with the model's formulation \cite{aliev1996}.
Finally, the myofiber architecture is prescribed using the atrial \ac{LDRBM} described in~\cite{piersanti2021modeling}, as depicted in \cref{fig:fiber}(b), by specifying \texttt{Geometry type = Left atrium} in the subsection \texttt{Fiber generation}. It is important to remark that in the \ac{APf} model formulation, the transmembrane potential $u$ is a dimensionless variable ranging from 0 to 1. However, for visualization purposes, the actual transmembrane potential $u_{\text{mV}}$ in millivolts is obtained by postprocessing the numerical results using the formula $u_{\text{mV}} = (100u - 80)~$[mV] \cite{aliev1996}. This conversion is applied to the data shown in Figure~\ref{fig:u}(b).

\subsubsection{Idealized left ventricle}
\label{subsec:ideal_lv}
We simulate the electrophysiology wavefront propagation in an idealized left ventricular hexahedral mesh, shown in Figure \ref{fig:mesh}(c), using the \acf{BO} ionic model \cite{buenoorovio2008}. The simulation can be run either by modifying the generic parameter file \texttt{lifex\_electrophysiology.prm} or by utilizing the dedicated parameter file \texttt{params\_ideal\_lv.prm}.
\begin{lstlisting}[language=prm]
$ ./lifex_electrophysiology-1.5.0-x86_64.AppImage \ 
      -f params_ideal_lv.prm
\end{lstlisting}
This test case can be executed at various levels of hierarchical grid refinements by modifying the parameter value \texttt{Number of refinements} in the \texttt{Mesh and space discretization} subsection. We remark that this is possible for all \lifex{} simulations conducted with hexahedral meshes, but is not available for tetrahedral meshes, due to lack of support in \dealii{}.
\begin{lstlisting}[language=prm]
subsection Electrophysiology
  subsection Mesh and space discretization
    set Element type          = Hex
    set Number of refinements = 3
    # ...
  end
  # ...
end
\end{lstlisting}
For this simulation, we use a second-order \ac{BDF} temporal scheme with a time step of $\Delta t=5 \times 10^{-5}$ s and a final time of $T=0.15$ s. The ionic model is initialized using 1000-cycle long single-cell simulations with a cardiac period of $0.8$ s. We employ a pacing protocol where three ventricular endocardial areas are activated using Gaussian impulses \cite{piersanti20223d0d}. The characteristics of these impulses, such as amplitude, duration, and initial time, can be set in the file parameter subsection \texttt{Applied current/Gaussian}.
\begin{lstlisting}[language=prm]
subsection Electrophysiology
  # ...
  subsection Applied current
    subsection Gaussian
      set Active         = true
      set Impulse sites  = -0.0271565 0.00506014 0.0141453, \
                           -0.0068242 -0.0187902 0.0382122, \ 
                            0.02695 0.00195906 0.0177283
      set Impulse amplitudes          = 300, 300, 300
      set Impulse standard deviations = 2.5e-3, 2.5e-3, 2.5e-3
      set Impulse initial times       = 0, 0, 0
      set Impulse durations           = 2e-3, 2e-3, 2e-3
      set Impulse period              = 0.8
    end
  end
# ...
end
\end{lstlisting}
The myofiber architecture is prescribed using the \ac{RL} ventricular \ac{LDRBM} \cite{piersanti2021modeling}, see \cref{fig:fiber}(c), by specifying \texttt{Geometry type = Left ventricle}  and setting \texttt{Algorithm type = RL} in the subsection \texttt{Fiber generation}.

The simulation results are reported in \cref{fig:u}(c) and \cref{fig:at}(c), where a snapshot of the transmembrane potential and the total activation time are displayed, respectively. Note that the transmembrane potential $u_{\text{mV}}$ shown in \cref{fig:u}(c) is obtained by postprocessing the numerical results using the formula $u_{\text{mV}}=(85.7u-84)~[mV]$~\cite{buenoorovio2008}.

\subsubsection{Ventricular slab}
\label{subsec:slab_lv}
This test case serves as an explanatory example for the multi-volume simulation framework of \lifexep{}. We simulate a portion of an idealized left ventricular geometry, also referred to as ventricular slab \cite{lifex-fiber}, where the computational domain is divided into three volumetric regions: Sub-endocardial, Myocardial, and Sub-epicardial layers, as shown in \cref{fig:mesh}(d).

To perform multi-volume simulations, the parameter template is created by specifying the volume labels on the command line as follows:
\begin{lstlisting}[language=prm]
$ ./lifex_electrophysiology-1.5.0-x86_64.AppImage -g \
      -vol "Sub Endocardium" "Myocardium" "Sub Epicardium"
\end{lstlisting}
Doing so, the file \texttt{lifex\_electrophysiology.prm} will contain three subsections named \texttt{Sub Endocardium}, \texttt{Myocardium} and \texttt{Sub Epiucardium}, where the volumetric tags (\texttt{Material IDs}) and all volume-specific parameters are prescribed. We use the \ac{TTP06} ionic model with a different \texttt{Cell type} for each layer (\texttt{Endocardium}, \texttt{Myocardium} and \texttt{Epicardium}) \cite{tentusscher2006}.
\begin{lstlisting}[language=prm]
subsection Electrophysiology
  # ...
  subsection Physical constants and models
    # ...
    subsection Sub Endocardium
      set Material IDs   = 1
      set Ionic model    = TTP06
      # ...
      subsection Ionic model parameters
        subsection TTP06
          set Cell type  = Endocardium
          # ...
        end
        # ...
      end
    end
    
    subsection Myocardium
      set Material IDs   = 2
      set Ionic model    = TTP06
      # ...
      subsection Ionic model parameters
        subsection TTP06
          set Cell type  = Myocardium
          # ...
        end
        # ...
      end
    end
    
    subsection Sub Epicardium
      set Material IDs   = 3
      set Ionic model    = TTP06
      # ...
      subsection Ionic model parameters
        subsection TTP06
          set Cell type  = Epicardium
          # ...
        end
        # ...
      end
    end
    # ...
  end
  # ...
end
\end{lstlisting}
The simulation can then be run by providing on the command line the same volume labels used when generating the parameter file:
\begin{lstlisting}[language=prm]
$ ./lifex_electrophysiology-1.5.0-x86_64.AppImage \ 
      -f params_slab_lv.prm \
      -vol "Sub Endocardium" "Myocardium" "Sub Epicardium"
\end{lstlisting}
A second-order \ac{BDF} temporal scheme is employed in this simulation, with a time step of $\Delta t=5 \times 10^{-5}$ s and a final time of $T=0.12$ s. The ionic model is initialized in each volume by conducting 1000-cycle long single-cell simulations with a cardiac period of $0.8$ s. One single-cell simulation is run for every volume, so that the initialization is consistent with the spatially-varying parameters. The simulation onset is obtained with a single spherical impulse, of radius $r=2.5 \times 10^{-3}$ m. The myofiber architecture is prescribed in the \texttt{Fiber generation} subsection by utilizing the \texttt{Slab} \ac{LDRBM}, as shown in \cref{fig:fiber}(d). The fiber orientations are set to exhibit a linear transmural variation from $60^{\circ}$ to $-60^{\circ}$, passing from the endocardial to the epicardial surface.

The simulation results are shown in \cref{fig:u}(d) and \cref{fig:at}(d), which display a snapshot of the transmembrane potential and the total activation time, respectively.

\subsection{Pathological electrophysiology simulations}
\label{subsec:phato}
In the following section, we present the pathological electrophysiology simulations applied to realistic left atrial and left ventricular geometries.

\subsubsection{Realistic pathological left atrium}
\label{subsec:real_la}
We perform a simulation of reentrant drivers typical of \ac{AF} in a realistic left atrial tetrahedral mesh \cite{courtemanche1998,roney2022predicting}, as shown in Figure \ref{fig:mesh}(e, h).

In this pathological scenario, four volumes have been introduced, named \texttt{Healthy}, \texttt{Fibrosis Mild}, \texttt{Fibrosis} and \texttt{Scar}. We use the \ac{CRN} ionic model adopting different ionic conductances and conductivity values to model the pathophysiological behaviour in the fibrotic regions, that are labelled according to their \ac{IIR} \cite{roney2022predicting}. Specifically, in the \texttt{Fibrosis Mild} volumetric region, we changed, with respect to the default \ac{CRN} values, the parameters related to the transient outward current conductance \texttt{gto}, the L-type intracellular $\text{Ca}^{2+}$ current conductance \texttt{gCaL}, and the delayed rectifier current represented by \texttt{gKur\_fix} and \texttt{gKur\_var}, to model the effects of chronic \ac{AF}. In the \texttt{Fibrosis} region, we set specific current conductance values for \texttt{gK1}, \texttt{gCaL} and \texttt{gNa}, which stand for the inward rectifier, the L-type calcium, potassium and sodium current conductances, respectively. Conductances and conductivities were adjusted according to \cite{zahid2016patient}. We also consider the \texttt{Healthy} physiological region, and the purely non-conductive region \texttt{Scar}, for which we set \texttt{Disable conduction = true}.
\begin{lstlisting}[language=prm]
subsection Electrophysiology
  # ...
  subsection Physical constants and models
    Time 0D simulation for variables initialization = 12
    subsection Healthy
      set Material IDs   = 1
      set Ionic model    = CRN
      subsection Monodomain conductivities
        set Longitudinal conductivity = 6.00e-4
        set Transversal conductivity  = 0.50e-4
        set Normal conductivity       = 0.50e-4
      end
    end
    
    subsection Fibrosis Mild
      set Material IDs   = 2
      set Ionic model    = CRN
      subsection Monodomain conductivities
        set Longitudinal conductivity = 2.000e-4
        set Transversal conductivity  = 0.175e-4
        set Normal conductivity       = 0.175e-4
      end
      # ...
      subsection Ionic model parameters
        subsection CRN
          subsection Physical constants
            #..
            set gto      = 0.0826
            set gCaL     = 0.037125
            set gKur_fix = 0.0025
            set gKur_var = 0.025
          end
          #..
        end
      end
    
      subsection Fibrosis
        set Material IDs   = 3
        set Ionic model    = CRN
        subsection Monodomain conductivities
          set Longitudinal conductivity = 0.500e-4
          set Transversal conductivity  = 0.050e-4
          set Normal conductivity       = 0.050e-4
        end
        # ...
        subsection Ionic model parameters
          subsection CRN
            subsection Physical constants
              #..
              set gNa      = 4.68
              set gK1      = 0.045
              set gCaL     = 0.061875
            end
            #..
          end
          # ...
        end
      end
    
      subsection Scar
        set Material IDs       = 4
        set Disable conduction = true
      end
      # ...
    end
    # ...
  end
  # ...
end
\end{lstlisting}
For the spatial discretization we employ second-order \acp{FE} by setting \texttt{FE space degree = 2}. We use a first-order \ac{BDF} temporal scheme with a time step $\Delta t=5~\times~10^{-5}$~s and a final time $T=1.5$ s. We initialize the ionic model by running a 24-cycle long single-cell simulation using an impulse period of $0.5$ s \cite{zahid2016patient}. We use a pacing protocol with a sequence of four multiple spherical impulses delivered every $220 \times 10^{-3}$ s. This timing is set under the parameter \texttt{set Impulse initial times} in the \texttt{Applied current/Spherical} subsection of the parameter file. Finally, the myofiber architecture is prescribed using the atrial \ac{LDRBM} \cite{piersanti2021modeling}, as depicted in  \cref{fig:fiber}(e), by specifying \texttt{Geometry type = Left atrium} in the subsection \texttt{Fiber generation}.

The simulation can be run using
\begin{lstlisting}[language=prm]
$ ./lifex_electrophysiology-1.5.0-x86_64.AppImage \ 
      -f params_real_la.prm \
      -vol "Healthy" "Fibrosis Mild" "Fibrosis" "Scar"
\end{lstlisting}
The simulation results are presented in \cref{fig:u}(e) and \cref{fig:at}(e), which display a snapshot of the transmembrane potential and the total activation time, respectively. In this video\footnote{\url{https://polimi365-my.sharepoint.com/:v:/g/personal/10594253_polimi_it/EcRfIyP2ya5EsD6IixPi_4ABCaZfJViTwwgeDpKBzCsquw?e=VUMMP9}} (online version) we show the evolution of the transmembrane potential.

\subsubsection{Realistic pathological left ventricle}
\label{subsec:real_lv}
We simulate a \acf{VT} macro-reentrant circuit in a realistic left ventricle \cite{mendoncacosta2021}, as shown in Figure \ref{fig:mesh}(f, i).

To account for the presence of grey zone and scar, we consider three volumes: \texttt{Healthy}, \texttt{Border Zone} and \texttt{Scar}. Different ionic conductances in the \ac{TTP06} ionic model and conductivity values in the monodomain equation are used to characterize the electrophysiology properties within each region. In particular, in the \texttt{Border Zone}, we reduce the conductances of the peak sodium, L-type calcium, and potassium currents  \texttt{GNa}, \texttt{GCaL}, \texttt{Gkr} and \texttt{Gks\_myo}, according to \cite{arevalo2016arrhythmia}. The \texttt{Healthy} region has a physiological configuration, while the \texttt{Scar} region is modeled as non-conductive (setting \texttt{Disable conduction = true}).
\begin{lstlisting}[language=prm]
subsection Electrophysiology
  # ...
  subsection Physical constants and models
    Time 0D simulation for variables initialization = 800
    subsection Healthy
      set Material IDs   = 1
      set Ionic model    = TTP06
      subsection Monodomain conductivities
        set Longitudinal conductivity = 9.0e-5
        set Transversal conductivity  = 1.8e-5
        set Normal conductivity       = 1.8e-5
      end
      # ...
    end
    
    subsection Border Zone
      set Material IDs   = 3
      set Ionic model    = TTP06
      subsection Monodomain conductivities
        set Longitudinal conductivity = 1.5e-5
        set Transversal conductivity  = 1.5e-5
        set Normal conductivity       = 1.2e-5
      end
      # ...
      subsection Ionic model parameters
        #..
        subsection TTP06
          set Cell type              = Myocardium
          subsection Physical constants
            #..
            set Gkr                  = 0.0459
            set Gks_myo              = 0.0294
            set GNa                  = 5.63844
            set GCaL                 = 0.00001234
          end
          #..
        end
      end
    end
    
    subsection Scar
      set Material IDs       = 4
      set Disable conduction = true
    end
    # ...
  end
  # ...
end
\end{lstlisting}
We use second-order \acp{FE} for the spatial discretization, with one level of hierarchical grid refinement (\texttt{Number of refinements = 1}), and a second-order \ac{BDF} temporal scheme with a time step of $\Delta t=5 \times 10^{-5}$ s and a final time of $T=2.5$ s. The ionic model is initialized using 1000-cycle long single-cell simulations with a cardiac period of $0.8$ s. We employ a pacing protocol with a sequence of multiple spherical impulses delivered in a specific ventricular endocardial area. The myofiber architecture is prescribed using the \ac{RL} ventricular algorithm \cite{piersanti2021modeling}, see \cref{fig:fiber}(f), by specifying \texttt{Geometry type = Left ventricle}  and setting \texttt{Algorithm type = RL} in the subsection \texttt{Fiber generation}.

The simulation can be run using
\begin{lstlisting}[language=prm]
$ ./lifex_electrophysiology-1.5.0-x86_64.AppImage \ 
      -f params_real_lv.prm \
      -vol "Healthy" "Border Zone" "Scar"
\end{lstlisting}
The simulation results are presented in \cref{fig:u}(f) and \cref{fig:at}(f), which display a snapshot of the transmembrane potential and the total activation time, respectively.  In this video\footnote{\url{https://polimi365-my.sharepoint.com/:v:/g/personal/10594253_polimi_it/EUOytZLaj1VGnsxoZxuZcb4BXSqKaGG63oUZPN6EPVt67g?e=77zmAP}} (online version) we show the evolution of the transmembrane potential.

\subsection{Strong scalability test}\label{subsec:scalability}
We consider the slab benchmark of \cite{niederer2011verification}, and discretize the domain with a structured hexahedral mesh of \num{47185920} elements and \num{47744577} nodes, corresponding to a mesh size of approximately $dx = \SI{3.6e-5}{\milli\meter}$. We set $\Delta t = \SI{1e-4}{\second}$ and $T = \SI{5e-2}{\second}$, and we use linear \acp{FE} and the second-order \ac{BDF} scheme for time discretization.
\begin{figure}[t!]
    \centering
    \includegraphics{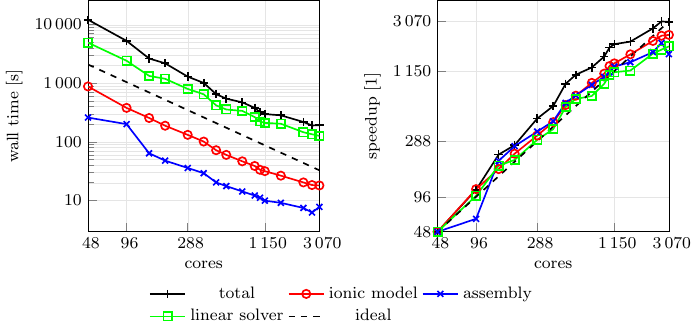}
    \caption{Computational time (left) and parallel speedup (right) against the number of cores for the strong scalability test. Dashed lines indicate the ideal linear scaling.}
    \label{fig:scalability}
\end{figure}
We run a strong scalability test, varying the number of parallel processes used in the computation and measuring the wall time necessary for the solution of the ionic models, the assembly of the monodomain system and the solution of the linear system. The test was run on the CINECA GALILEO100 supercomputer.

The wall times, plotted in \cref{fig:scalability} for the different steps of the solver, scales linearly up to over \num{1000} cores, confirming the results of \cite{africa2022lifex} on the scalability properties of the \lifex{} core components. We report in \cref{tab:performance} a breakdown of the computational cost of the different sections for the simulation with \num{960} cores. Most of the computational time is spent for the solution of the linear system \eqref{eqn: linear-system}. The matrix of the linear system is assembled only once (since it is the same at every time step), and the right-hand-side is efficiently recomputed at every time step by means of matrix-vector products, resulting in a very small computational cost for the system assembly phase. Moreover, since the ionic model is solved independently at each degree of freedom, the associated computational cost scales almost perfectly, and becomes very small if a sufficiently large number of parallel processes is employed.

\begin{table}[]
    \centering
    \scalebox{0.8}{
        \begin{tabular}{c S S S}
            Section & {Wall time [\si{\second}]} & {Wall time per time step [\si{\second}]} & {Relative wall time [\si{\percent}]} \\
            \hline
            Total               & 379.51 & {-}  & 100.0 \\
            Linear solver       & 226.98 & 0.45 & 71.0 \\
            Ionic model update  &  45.13 & 0.09 & 14.1 \\
            Initialization      &  33.76 & {-}  & 10.6 \\
            Monodomain assembly &  10.33 & 0.02 & 3.2 \\
            Other	            &   3.55 & 0.01 & 1.1 \\
            \hline
    \end{tabular}}
    \caption{Summary of the computational costs for the strong scalability test, using \num{960} parallel cores. For each section, we report the total wall time, the wall time for each time step and the wall time relative to the total. Sections are sorted in descending order of cost.}
    \label{tab:performance}
\end{table}

\subsection{Output and visualization}
\label{subsec:output}
Two different types of output file formats are available in the \lifexep{} release: \texttt{HDF5} and \texttt{csv}. Both of them can be enabled and configured in the \texttt{Output} subsections:
\begin{lstlisting}[language=prm]
subsection Electrophysiology
  # ...
  subsection Output
    set Enable output     = true
    set Filename          = solution
    set Enable CSV output = true
    set CSV filename      = electrophysiology.csv
  end
  # ...
end
\end{lstlisting}
The \texttt{HDF5} output is available in the following subsections:
\begin{itemize}
    \item \texttt{Electrophysiology/Output},
    \item \texttt{Electrophysiology/Activation time},
    \item \texttt{Fiber generation/Output}.
\end{itemize}
This generates an XDMF file named \verb|output_filename.xdmf| (which links to a corresponding HDF5 output file \verb|output_filename.h5|). These files can be visualized using \texttt{ParaView}\footnote{\url{https://www.paraview.org}}, an open-source multi-platform data analysis and visualization application, see e.g., \cref{fig:fiber}, \cref{fig:u} and \cref{fig:at}. The \texttt{HDF5} format ensures that the output can be easily post-processed, not only for visualization purposes but also as input for more advanced computational pipelines.

The Comma-Separated Values (\texttt{csv}) format consists of delimited text files where values are separated by commas, with each line representing a specific data record. The \texttt{csv} files can be found in different subsections: 
\begin{itemize}
    \item \texttt{Ionic model parameters/Output},
    \item \texttt{Ionic model parameters/0D Output},
    \item \texttt{Electrophysiology/Output}.
\end{itemize}
These \texttt{csv} output files can be conveniently used to plot electrophysiology variable (min, max and pointwise) values over time in the computed numerical simulation.

\section{Conclusions}\label{sec:conclusions}

In this work, we introduced \lifexep{}, a robust and advanced software specifically designed for simulating the electrophysiology activity of the cardiac muscle. With the goal of addressing the computational challenges associated with cardiac simulations, \lifexep{} provides efficient numerical methods while maintaining precision and accuracy. \lifexep{} incorporates a numerical solver for the monodomain equation coupled with both phenomenological and second-generation ionic models, namely \acl{APf} \cite{aliev1996}, Bueno-Orovio \cite{buenoorovio2008}, \ac{TTP06} \cite{tentusscher2006}, and \ac{CRN} \cite{courtemanche1998}. These models are discretized in time using the \ac{BDF} scheme and in space using the \ac{FE} method of orders $1$ and $2$ on tetrahedral meshes, and of arbitrary degree on hexahedral meshes, thus providing a comprehensive framework for modeling the electrical activity of the heart under both physiological and pathological conditions.

Leveraging the capabilities of \lifex{}, \lifexep{} provides users with a user-friendly and flexible interface, facilitated by self-documenting parameter files for easy simulation setup. For enhanced accessibility, \lifexep{} is distributed in an \texttt{AppImage} binary format, rendering it universally compatible with any recent \texttt{x86-64 Linux} system. 
Researchers from diverse backgrounds, such as medicine and bio-engineering, can readily access and utilize \lifexep{} for in-silico simulations. The underlying principles and structure of \lifexep{} can be readily understood thanks to the comprehensive technical and mathematical documentation.

As unique and distinctive features, \lifexep{} provides two options for prescribing myocardial fibers. Users can either import them from a file or generate them online by exploiting the \acp{LDRBM} presented in \cite{piersanti2021modeling} and implemented in the previous release \lifexfiber{} \cite{lifex-fiber}. Moreover, it supports spatial heterogeneity in the choice of both models and physical coefficients, easily configurable through a convenient parameter file, without the need to access and modify the source code.

\lifexep{} benefits from its high-performance computing capabilities, achieving ideal parallel speedup on thousands of cores. The accuracy and reliability of \lifexep{} have been verified through a benchmark for computational electrophysiology, ensuring the validity of its results. Furthermore, a range of idealized and realistic cardiac simulations in both physiological and pathological settings highlights its capabilities and versatility in capturing complex cardiac dynamics and its potential for patient-specific studies. \lifexep{} offers the capability to facilitate the simulation of pathological scenarios, allowing the creation of scars, grey zones, and an arbitrary number of conduction "levels". This potential impact in the study of pathologies can prove to be highly valuable.

In conclusion, \lifexep{} 
provides to the scientific community a comprehensive, high-performance, and user-friendly software for conducting in-silico cardiac electrophysiology simulations.

In the future, efforts will be made to further improving the accuracy and efficiency of \lifexep{}. One possible approach is the adoption of lookup tables instead of repeatedly evaluating the expensive functions in the ionic current term, which could provide a speed-up, although with a potential impact on accuracy~\cite{cooper2015cellular}. Inexact solvers, such as those based on domain decomposition methods, show promising optimality and scalability properties \cite{zampini2014inexact}. Recent studies have also highlighted the advantages of employing high-order discretization schemes for accurately capturing the intricate electrical wavefront propagation observed in cardiac electrophysiology~ \cite{africa2023matrix}. These schemes not only offer improved accuracy but also enable the implementation of efficient matrix-free solvers that require minimal memory usage~\cite{africa2023matrix}. This opens up the possibility of leveraging GPU architectures for accelerated computations \cite{del2022fast}. Additionally, the use of higher-order or adaptive time stepping schemes and \(hp\)-adaptive \ac{FE} holds the potential to achieve greater accuracy in simulations while optimizing computational efficiency \cite{CHAMAKURI2022295}. These developments aim to further enhance the capabilities of \lifexep{} and expand its range of features.

\section*{List of abbreviations}
\begin{acronym}
    \acro{FE}{Finite Element}
    \acro{BDF}{Backward Differentiation Formula}
    \acro{TTP06}{ten Tusscher-Panfilov 2006}
    \acro{APf}{Aliev-Panfilov}
    \acro{BO}{Bueno-Orovio}
    \acro{CRN}{Courtemanche-Ramirez-Nattel}
    \acro{AF}{Atrial Fibrillation}
    \acro{VT}{Ventricular Tachycardia}
    \acro{CRT}{Cardiac Resynchronization Therapy}
    \acro{MRI}{Magnetic Resonance Imaging}
    \acro{LGE-MRI}{Late Gadolinium Enhancement-Magnetic Resonance Imaging}
    \acro{AP}{Action Potential}
    \acro{CV}{Conduction Velocity}
    \acro{AF}{Atrial Fibrillation}
    \acro{VT}{Ventricular Tachycardia}
    \acro{LBBB}{Left Bundle Branch Block}
    \acro{ECG}{Electrocardiogram}
    \acro{BSPM}{Body Surface Potential Mapping}
    \acro{CT}{Computed Tomography}
    \acro{ODEs}{Ordinary Differential Equations}
    \acro{PDE}{Partial Differential Equation}
    \acro{ICI}{ionic current interpolation}
    \acro{RL}{Rossi-Lassila}
    \acro{LDRBM}{Laplace-Dirichlet Rule-Based Method}
    \acro{IIR}{Imaging Itensity Ratio}
    \acro{MPI}{Message Passing Interface}
\end{acronym}

\section*{Availability and requirements}
\textbf{Project name}: \lifexep{} \\
\textbf{Project home page}: \url{https://lifex.gitlab.io/} \\
\textbf{Operating system(s)}: \texttt{Linux} (\texttt{x86-64})\\
\textbf{Programming language}: \texttt{C++} \\
\textbf{Other requirements}: \texttt{glibc} version \texttt{2.28} or higher \\
\textbf{License}: CC BY-NC-ND 4.0 \\
\textbf{Any restrictions to use by non-academics}: no additional restriction.

\section*{Availability of data and materials}
All input data, meshes and the binary executable of \lifexep{} can be found at \zenodourl{}.

\section*{Acknowledgements}\label{sec:funding}
This project has received funding from the European Research Council (ERC) under the European Union's Horizon 2020 research and innovation program (grant agreement No 740132, iHEART - An Integrated Heart Model for the simulation of the cardiac function, P.I. Prof. A. Quarteroni). P.C.A, R.P., F.R., M.B., M.F., S.P. and L.D. are members of the INdAM research group GNCS. 

\section*{Authors' contributions}
\textbf{P.C.A.}: conceptualization, methodology, software (development and maintenance), formal analysis, supervision, writing (original draft).
\textbf{R.P.}: conceptualization, methodology, software (development and simulation), formal analysis, writing (original draft).
\textbf{F.R.}: conceptualization, methodology, software (development and simulation), formal analysis, writing (original draft).
\textbf{M.B.}: conceptualization, methodology, software (development, simulation, testing and maintenance), formal analysis, writing (original draft). 
\textbf{M.S.}: conceptualization, methodology, software (development and simulation), formal analysis, writing (original draft).
\textbf{M.F.}: conceptualization, methodology, software (development), formal analysis, writing (review and editing).
\textbf{S.P.}: conceptualization, methodology, software (development and simulation), formal analysis, writing (original draft).
\textbf{L.D.}: project administration, writing (review and editing), supervision.
\textbf{A.Q.}: funding acquisition, project administration, writing (review and editing), supervision.

\printbibliography

\end{document}